\documentclass[lettersize,journal]{IEEEtran}
\usepackage{amsmath,amsfonts}
\usepackage{algorithm}
\usepackage{algpseudocode} 
\usepackage{array}
\usepackage[caption=false,font=normalsize,labelfont=sf,textfont=sf]{subfig}
\usepackage{textcomp}
\usepackage{stfloats}
\usepackage{url}
\usepackage{verbatim}
\usepackage{graphicx}
\usepackage{cite}
\usepackage{mathtools}
\usepackage{booktabs}
\usepackage{multirow}
\usepackage{xcolor}
\usepackage{lipsum}
\usepackage[inline]{enumitem}
\usepackage{amssymb}
\usepackage{cuted}
\usepackage{siunitx}
\usepackage{amsthm} 
\newtheoremstyle{customremark}
  {\topsep}        
  {\topsep}        
  {}               
  {}               
  {\bfseries}      
  {.}              
  {.5em}           
  {}               

\theoremstyle{customremark}
\newtheorem{remark}{Remark}

\algrenewcommand\algorithmicrequire{\textbf{Input:}}
\algrenewcommand\algorithmicensure{\textbf{Output:}}

\algnewcommand{\Inputs}[1]{%
  \State \textbf{Inputs:}
  \Statex \hspace*{\algorithmicindent}\parbox[t]{.85\linewidth}{\raggedright #1}
}

\algnewcommand{\Initialize}[1]{%
  \State \textbf{initialization:}~#1
}

\algnewcommand{\Do}{\State \textbf{do}}

\algnewcommand{\DoWhile}[1]{\State \textbf{while}~#1}

\usepackage[capitalize]{cleveref}
\crefname{section}{Sec.}{Secs.}
\Crefname{section}{Section}{Sections}
\Crefname{table}{Table}{Tables}
\crefname{table}{Tab.}{Tabs.}

\usepackage{titlesec}
\titleformat{\subsubsection}[runin]     
  {\normalfont}                         
  {\thesubsubsection)}                  
  {0.5em}                               
  {\bfseries}[:]                        
\titlespacing*{\subsubsection}{0pt}{0.5em}{0.5em}  
\renewcommand\thesubsubsection{\arabic{subsubsection}} 

\DeclareMathOperator*{\argmin}{arg\,min} 

\allowdisplaybreaks

\begin{document}

\title{Efficient sparse GP-MPC with accurate mean and variance propagation applied for quadcopter flight control}

\author{Giannis Badakis$^{1}$,~
Mircea Lazar$^{1}$,~and
Roland T{\'o}th$^{1,2}$
\thanks{$^{1}$Giannis Badakis, Mircea Lazar, and Roland T{\'o}th are with the Control Systems Group, Eindhoven University of Technology, Eindhoven, the Netherlands (e-mail: johnbdks@gmail.com; m.lazar@tue.nl; r.toth@tue.nl).}
\thanks{$^{2}$Roland T{\'o}th is also with the Systems and Control Laboratory, HUN-REN Institute for Computer Science and Control, Budapest, Hungary, and with the Vehicle Industry Research Center, Sz{\'e}chenyi Istv{\'a}n University, Gy\H{o}r, Hungary.}
}

\maketitle

\begin{abstract}
This paper presents a computationally efficient approach for Gaussian process model predictive control (GP-MPC), where Gaussian process (GP) regression is used to complement a baseline model of the system dynamics. The proposed method achieves propagation of both the predicted mean and variance, thereby significantly reducing conservativeness compared with existing GP-MPC formulations. The nonlinear GP-MPC problem is reformulated into an exact linear parameter-varying (LPV) structure that preserves the nonlinear prediction dynamics in affine form without introducing further approximation. Moreover, closed-form derivations of moment matching (MM) predictions for sparse GPs are developed, including both mean and variance propagation under uncertain inputs, which improves scalability to larger datasets. This further enables recasting the resulting GP-MPC problem as a sequence of quadratic programs (QPs), which can be solved efficiently. The proposed framework significantly improves runtime efficiency while maintaining prediction accuracy, as demonstrated through simulation and real-world experiments on a Crazyflie 2.1 micro quadcopter.

\end{abstract}


\section{Introduction}
\label{sec:introduction}

\IEEEPARstart{M}{odel} predictive control (MPC) has become a cornerstone of modern control applications, providing a unified framework for optimizing performance over a receding horizon while systematically handling input and state constraints. However, its effectiveness fundamentally depends on the fidelity of the prediction model. In situations where first-principles models are inaccurate, difficult to derive, or computationally expensive to evaluate, learning-based models have gained increasing attention as a viable alternative~\cite{RobotLearning}.

A widely adopted and effective strategy for model learning is model augmentation~\cite{BackFlipping, HybridModelling, SaKaAr23_AgileMavsANN_MPC}, where first-principles models are enhanced with learned corrections. Typical machine learning approaches used to capture unmodeled phenomena include~\emph{artificial neural networks} (ANNs)\cite{BeScTo23_SUBSPACE, GaAlDa23_AnnSurvey} and probabilistic~\emph{Gaussian processes} (GPs)~\cite{CaCaLi23_KernelMethods}. An advantage of GPs over ANNs is their ability to provide predictive variance in addition to a mean prediction, which quantifies model uncertainty.
This property is particularly critical in control applications where uncertainty quantification is essential for safe decision-making and accurate closed-loop control performance is required~\cite{LazarGuarantees}.

The combination of GP models with the MPC framework, commonly referred to as GP-MPC~\cite{HeKaZe20_GPMPC1}, enables principled control synthesis that accounts for model uncertainty during planning. This is achieved by using GP predictions to forecast future states while incorporating predictions of their uncertainty into the cost and constraint formulation.

However, a key challenge in GP-MPC lies in the recursive inference through the augmented model to predict future states within the prediction horizon. As predictions extend further into the future, the complexity of the state probability distributions increases, resulting in an intractable stochastic optimization problem typically addressed through a series of approximations. For this purpose, two common approaches are: 
\begin{enumerate*}[label=(\roman*)]
\item applying a Taylor series expansion to the predictive GP mean and covariance functions~\cite{GiRaMu02_Uncertain_GP_Taylor}, or
\item computing the first two moments analytically---known as~\emph{moment matching} (MM)~\cite{QuGiRa03_Uncertain_GP_EMM, DeHuHa_AnalyticMM, De_PILCO_PhD}---under the assumption of Gaussian inputs and a \emph{squared exponential} (SE) kernel\cite{De_PILCO_PhD}.
\end{enumerate*}
Both methods enable the formulation of a deterministic~\emph{nonlinear} MPC (NMPC) problem that propagates only the mean and variance dynamics and can be solved using standard nonlinear programming techniques.

In practice, most GP-MPC implementations~\cite{HeLiZe18_GPMPC2, CaArWe19_MPCArm, ToKaFo21_GPMPCDrone} adopt the Taylor series approximation to enable real-time policy generation. While this reduces computational complexity, it introduces certain limitations. In particular, Taylor-based methods neglect variance-dependent terms in the GP predictive mean, which may lead to conservative control policies and reduced performance. The approximation error may accumulate over the prediction horizon, affecting prediction accuracy. Moreover, for moderate- or high-order systems such as those considered in the aforementioned works, the resulting NMPC problem is often computationally too demanding for real-time operation. As a result, many GP-MPC implementations omit covariance propagation altogether, which can reduce robustness and lead to suboptimal use of the GP model. Although pragmatic strategies such as fixing the covariance dynamics across the prediction horizon using the previous optimal trajectory~\cite{HeLiZe18_GPMPC2, AmAnAn23_ZOH} embed variance information in the constraints, they still rely on Taylor expansion for propagating the model within the MPC framework and typically neglect the effect of changing variance on the evolution of the mean dynamics.

An alternative approach was introduced in~\cite{PoPeTo23_LPVGPMPC}, where a~\emph{linear parameter-varying} (LPV) formulation significantly reduced the computational burden of NMPC by solving a sequence of~\emph{quadratic programs} (QPs). To the best of the authors’ knowledge, this formulation enabled, for the first time in a GP-MPC framework, the efficient use of the MM approximation to more accurately approximate the GP predictive distribution and to embed variance terms in the mean dynamics. However, the LPV formulation in~\cite{PoPeTo23_LPVGPMPC} is obtained via a Jacobian-based first-order affine representation of the nonlinear prediction dynamics. Therefore, the LPV embedding is only locally exact around the scheduling trajectory, as higher-order nonlinear terms of the original nonlinear model are not preserved. In addition, the covariance dynamics were fixed with respect to the scheduling variables, potentially limiting the model adaptability to variations in the uncertainty due to the optimized input sequence. Furthermore, this LPV formulation includes all data points within the GP and the subsequent MM approximation, which restricts the scalability of the method to larger datasets. These factors limit the ability of the approach to adapt to complex and uncertain environments, eventually introducing conservatism in the control decisions.

In summary, although LPV formulations enable faster GP-MPC computation, Jacobian-based LPV embeddings of the nonlinear prediction dynamics are only locally exact and may therefore reduce model fidelity. They also fail to capture uncertainty propagation explicitly, since the LPV-formulation-based variance dynamics vanish as they are fully expressed by the scheduling trajectories. This work addresses both limitations through the following main contribution:

\begin{enumerate}[label = {(C\arabic*)}]
\item An approximation-free LPV-MPC formulation of the nonlinear GP-MPC problem under both Taylor series and MM propagation of the state distributions. \label{contributions:main}
\end{enumerate}

This formulation recasts the nonlinear GP prediction dynamics---both mean and variance---into an LPV structure without introducing additional approximation beyond the GP propagation method. In contrast to iterative schemes based on local linearization of the nonlinear optimization problem, the proposed reformulation preserves the original nonlinear prediction dynamics in an affine LPV form. By retaining model fidelity and the influence of uncertainty in the mean dynamics, it enables efficient and accurate iterative QP-based solutions to the GP-MPC problem. Moreover, the formulation is extended to sparse GPs, thereby improving scalability to larger datasets and broadening the applicability of the method.

Based on these, the following sub-contributions are made:
\begin{enumerate}[resume, label = {(C\arabic*)}]

\item Derivation of closed-form MM 
predictions (both mean and variance dynamics) for scalable sparse GPs under uncertain inputs. \label{contributions:item-one}

\item Adaptation of the iterative QP-based solution scheme of GP-MPC from~\cite{PoPeTo23_LPVGPMPC} to the proposed LPV-MPC formulation, where scheduling sequences are updated via the optimal state and input trajectories at each iteration.\label{contributions:item-two}

\item Efficient variance propagation within the LPV-based GP-MPC structure, enabling  uncertainty handling across the prediction horizon. \label{contributions:item-three}

\item Implementation of the developed LPV GP-MPC framework in real-time control of a Crazyflie 2.1 micro quadcopter. \label{contributions:item-four}

\end{enumerate}
Notably, the reformulated variance dynamics from sub-contribution~\ref{contributions:item-three} enable pre-computation of the uncertainty trajectory using prior QP solutions from sub-contribution~\ref{contributions:item-two}. This facilitates fast, sub-optimal solutions of the GP-MPC problem that remain tractable even for large-scale or high-order systems, while maintaining prediction accuracy without incurring additional computational cost compared with prior approaches such as~\cite{HeLiZe18_GPMPC2}.

The remainder of this paper is structured as follows:~\Cref{sec:preliminaries} presents the problem formulation and GP-based modeling.~\Cref{sec:mpc} reviews GP-MPC and introduces sparse MM, addressing~\ref{contributions:item-one}.~\Cref{sec:lpv-sol} develops the proposed LPV-based solution, highlighting the main contribution~\ref{contributions:main} and detailing~\ref{contributions:item-two} and~\ref{contributions:item-three}.~\Cref{sec:application} demonstrates the efficiency of the proposed scheme using both simulation results and real-world experiments on adaptive flight control of a Crazyflie 2.1 micro quadcopter, addressing~\ref{contributions:item-four}, and~\Cref{sec:conclusions} concludes the paper.

\section{Preliminaries \& problem statement}
\label{sec:preliminaries}
In this section, we introduce the notation used throughout the paper, formulate the considered system together with the predictive control and learning problems addressed in this work, and introduce model augmentation via GP regression and its sparse approximation based on inducing points.

\subsection{Notation}
\label{subsec:notation}
Let  $\mathbb{Z}$, $\mathbb{N}$, and $ \mathbb{R} $ denote the sets of integers, non-negative integers, and real numbers, respectively. Let $ \mathbb{S}^n $ be the set of symmetric matrices in $ \mathbb{R}^{n \times n} $, and define $ \mathbb{S}^n_{\succeq} = \{ X \in \mathbb{S}^n \mid X \succeq 0 \} $ and $ \mathbb{S}^n_{\succ} = \{ X \in \mathbb{S}^n \mid X \succ 0 \} $ as the sets of symmetric positive semidefinite and positive definite matrices, respectively. For a vector $ x \in \mathbb{R}^n $ and a matrix $ Q \in \mathbb{S}^n_{\succ} $, we define the $ Q $-weighted 2-norm as $ \lVert x \rVert_Q^2 = x^\intercal Q x $. Given matrices $X_1, \ldots, X_N$ of compatible dimensions, we define $\mathrm{col}(X_1, \ldots, X_N) := [X_1^\intercal \ \cdots \ X_N^\intercal]^\intercal$ as their vertical concatenation, and $\mathrm{diag}(X_1, \ldots, X_N)$ as their block-diagonal composition. For a matrix $X$, we denote by $\operatorname{vec}(X)$ the vector obtained by stacking its columns into a single column vector. Let $\mathbb{I}_{\tau_1}^{\tau_2} = \{ i\in\mathbb{Z} \mid \tau_1 \leq i \leq \tau_2 \}$ be a discrete index set. For a time index $k \in \mathbb{N}$, we denote the predicted value of $x(k+i)$ at time $k$ as $x(i|k)$. If $x$ and $y$ are random variables, we denote their expected value and variance by $\mathbb{E}\{x\}$ and $\mathrm{var}\{x\}$, respectively, and their covariance by $\mathrm{cov}\{x, y\}$. A Gaussian random variable $x$ with mean $\mu$ and variance $\Sigma$ is denoted as $x \sim \mathcal{N}(\mu, \Sigma)$.

\subsection{Considered system and problem setting}
\label{subsec:sys-model}
Consider the discrete-time nonlinear system defined as
\begin{equation}
    x(k+1) = f\bigl(x(k),u(k)\bigl) + g_\mathrm{d}\bigl(x(k),u(k)\bigr) + v(k),
    \label{eq:discrete_system}
\end{equation}
where $k \in \mathbb{N}$ is the discrete-time index, $x(k) \in \mathbb{R}^{n_\mathrm{x}}$ is the state, and $u(k) \in \mathbb{R}^{n_\mathrm{u}}$ is the input. The function $f : \mathbb{R}^{n_\mathrm{x}}\times\mathbb{R}^{n_\mathrm{u}} \mapsto \mathbb{R}^{n_\mathrm{x}}$ represents the discretized nominal model (e.g., a first-principles model of the system) and incorporates prior knowledge of the system dynamics, while the function $g_\mathrm{d} : \mathbb{R}^{n_\mathrm{x}}\times\mathbb{R}^{n_\mathrm{u}} \mapsto \mathbb{R}^{n_\mathrm{x}}$ represents unmodeled phenomena. Also, we assume that the external disturbance $v(k) \in \mathbb{R}^{n_\mathrm{x}}$ is an \emph{independent and identically distributed} (i.i.d.) \emph{white Gaussian noise} (WGN) process with $v(k)\sim \mathcal{N}\bigl(0_{n_\mathrm{x}}, \Sigma_v\bigr)$ where $\Sigma_v = \mathrm{diag}(\sigma_{v,1}^2,\dots,\sigma_{v,n_\mathrm{x}}^2)$. To simplify the notation of the arguments of both $f$ and $g_\mathrm{d}$, let $w := \mathrm{col} \bigl( x,u\bigr) \in \mathbb{R}^{n_\mathrm{w}}$ with $n_\mathrm{w} = n_\mathrm{x} + n_\mathrm{u}$.

We consider the problem of predictive reference tracking for the nonlinear system~\eqref{eq:discrete_system}. The objective is to design an MPC scheme that leverages an accurate model of the dynamics to efficiently track a desired trajectory over a finite horizon $N_{\mathrm{p}}$.

However, since only the baseline first-principles model $f$ is available, we formulate a second problem setting to identify the unknown function $g_\mathrm{d}$ from the residual dynamics $z(k) = x(k+1) - f\bigl(w(k)\bigr) = g_\mathrm{d}\bigl(w(k)\bigr) + v(k)$, with $z(k) \in \mathbb{R}^{n_\mathrm{z}}$ and $n_\mathrm{z} = n_\mathrm{x}$ based on a dataset $\mathcal{D}_N = \{(w_\tau, z_\tau)\}_{\tau=1}^N$ of $N$ data samples obtained from the data-generating system~\eqref{eq:discrete_system}.

The resulting augmented model must provide the MPC scheme with accurate and tractable multi-step predictions of the state evolution, i.e., $x(1|k), \dots, x(N_{\mathrm{p}}|k)$, along with quantified uncertainty, given the initial state $x(k) = x(0|k)$ and a sequence of future inputs $u(0|k), \dots, u(N_{\mathrm{p}}-1|k)$.

\subsection{Augmentation via Gaussian process learning}
\label{subsec:gp-learning}
We employ GP regression~\cite[Chap.~2]{RaCh05_GP} to address the modeling of $g_\mathrm{d}$ based on the observed residuals $z$ and the dataset $\mathcal{D}_N$ defined above. Specifically, we place a distribution over functions by modeling $g_\mathrm{d}$ with a GP surrogate, i.e., $g \sim \mathcal{GP}\left(m_\mathrm{d}, \kappa_\mathrm{d} \right)$, where $m_\mathrm{d}: \mathbb{R}^{n_\mathrm{w}} \mapsto \mathbb{R}^{n_\mathrm{z}}$ is a mean function encoding prior knowledge about $g_\mathrm{d}$, and $\kappa_\mathrm{d}: \mathbb{R}^{n_\mathrm{w}} \times \mathbb{R}^{n_\mathrm{w}} \mapsto \mathbb{R}^{n_\mathrm{z} \times n_\mathrm{z}}$ is a matrix-valued kernel function capturing smoothness assumptions and defining the function class.

Under this formulation, a vectorial $\mathcal{GP}: \mathbb{R}^{n_\mathrm{w}}\mapsto\mathbb{R}^{n_\mathrm{z}}$ defines a distribution over functions such that for any finite collection of inputs $\{w_\tau\}_{\tau=1}^N \subset \mathbb{R}^{n_\mathrm{w}}$, the associated outputs $\{g(w_\tau)\}_{\tau=1}^N$ follow a joint multivariate Gaussian distribution. We assume that each $g_\mathrm{i}(w)$ element of $g(w)$, with $i \in \mathbb{I}_{1}^{n_\mathrm{z}}$, is conditionally independent given $w$, so GP learning is simplified as independent scalar regressions $g = \mathrm{col}(g_1, \dots, g_{n_\mathrm{z}})$.

Accordingly, for each scalar prior $g_i \sim \mathcal{GP}(m_{\mathrm{d},i},\kappa_{\mathrm{d},i})$, we assume zero mean (i.e., $m_{\mathrm{d},i}=0$) and define an SE kernel:
\begin{equation}
    \kappa_{\mathrm{d},i}(w, \tilde w) = \sigma^2_i \exp \bigl\{ -\tfrac{1}{2} (w - \tilde w)^\intercal\Lambda^{-1}_i(w - \tilde w) \bigr\}, \label{eq:se_kernel}
\end{equation}
where $\sigma^2_i$ is the signal variance (also known as amplitude), and $\Lambda_i = \mathrm{diag}(\lambda_{i,1},\dots,\lambda_{i,n_\mathrm{w}}) \in \mathbb{S}^{n_\mathrm{w}}_{\succ}$ contains the length-scales along each input $w$ dimension. Since the training residuals $z_i$ include WGN process noise $v$, we incorporate this into the GP covariance structure by defining the total kernel $\kappa_i$ as:
\begin{equation}
\kappa_i(w, \tilde w) = \kappa_{\mathrm{d},i}(w, \tilde w) + \delta(w, \tilde w)\sigma_{v,i}^2, \label{eq:kernel}
\end{equation}
where $\delta$ is the Kronecker delta. Hyperparameters for each individual $g_i$ are collected in $\theta_i = \mathrm{col}\left(\sigma_i^2, \sigma_{v,i}^2, \lambda_{i,1}, \dots, \lambda_{i,n_\mathrm{w}}\right)$. Given the dataset $\mathcal{D}_{N,i} = \{(w_\tau, z_{i,\tau})\}_{\tau=1}^N$, with $Z_i = \mathrm{col}\left(z_{i,1},\dots,z_{i,N}\right)$ and $W = \mathrm{col}\left(w_1,\dots,w_N\right)$, the hyperparameters are estimated as
\begin{equation}
\theta_i^\star = \argmin_{\theta_i} \bigl(-\log p(Z_i \mid W, \theta_i)\bigr),
\end{equation}
where $-\log p(Z_i \mid W, \theta_i)$ denotes the negative log marginal likelihood, given by
\begin{equation}
-\log p(Z_i \mid W, \theta_i)
=
\frac{1}{2}
Z_i^\intercal K_{ww,i}^{-1}(\theta_i) Z_i
+
\frac{1}{2}\log \mathrm{det}K_{ww,i}(\theta_i).
\end{equation}

Adjusting the prior $g_i \sim \mathcal{GP}(0, \kappa_{\mathrm{d},i})$ based on the likelihood of the observed data $\mathcal{D}_{N,i}$, to predict a test output $z_i^\ast$ at a new query point $w^\ast$, we define a joint Gaussian distribution over $z_i^\ast$ and the training outputs $Z_i$ as 
\begin{equation}
    \begin{bmatrix}
        Z_i \\ z_i^\ast
    \end{bmatrix} \sim
    \mathcal{N}\Biggl(
    \begin{bmatrix}
    0_{N} \\ 0
    \end{bmatrix},
    \begin{bmatrix}
    K_{ww,i} & K_{w,i}(w^\ast) \\
    K_{w,i}^\intercal(w^\ast) & \kappa_i(w^\ast,w^\ast)
    \end{bmatrix}
    \Biggr),
\end{equation}
where $K_{ww,i}$ is the Gramian matrix in which $\left[ K_{ww,i} \right]_{\tau, \tilde{\tau}} = \kappa_i(w_\tau, w_{\tilde{\tau}})$, and $K_{w,i}(w^\ast)$ is a so-called kernel slice, with $\left[ K_{w,i}(w^\ast) \right]_{\tau} = \kappa_i(w_\tau, w^\ast)$ for $\tau, \tilde{\tau} \in \mathbb{I}_{1}^{N}$. Then, the conditional probability density function of $z_i^\ast$, that is $p(z^\ast_i|\mathcal{D}_{N,i}, w^\ast) = \mathcal{N}\bigl(\mu_{z,i}(w^\ast), \sigma_{z,i}^2(w^\ast)\bigr)$, is the predictive distribution of $z_i^\ast$ which is characterized by
\begin{subequations}
    \begin{align}
        \mu_{z,i}(w^\ast) &= K_{w,i}^\intercal(w^\ast) \overbrace{K_{ww,i}^{-1} Z_i}^{\alpha_i} \label{eq:pred_mean}, \\
        \begin{split}
            \sigma_{z,i}^2(w^\ast) &= \kappa_i(w^\ast,w^\ast) - K_{w,i}^\intercal(w^\ast) K_{ww,i}^{-1} K_{w,i}(w^\ast).
        \end{split}
        \label{eq:pred_var}
    \end{align}
    \label{eq:predictive_distr}
\end{subequations}

The complete multivariate predictive distribution of $z^\ast$ is given by the mean $\mu_z(w^\ast)$, where $ \mu_z = \mathrm{col}\left(\mu_{z,1},\dots,\mu_{z,n_\mathrm{z}}\right)$, and the variance $\Sigma_z(w^\ast)$, where $ \Sigma_z = \mathrm{diag}\left(\sigma_{z,1}^2,\dots, \sigma_{z,n_\mathrm{z}}^2\right)$.

Then, the discrete-time nonlinear system~\eqref{eq:discrete_system} can be expressed in the augmented form
\begin{equation}
    x(k+1) = f\bigl(w(k)\bigr) + \hat{z}(k),
    \label{eq:aug_model}
\end{equation}
where $\hat{z}(k) \sim \mathcal{N}\bigl(\mu_z(w(k)), \Sigma_z(w(k))\bigr)$ denotes the predictive distribution of the residual term $z(k)$ at the current $w(k)$, capturing both the unknown dynamics $g_\mathrm{d}(w(k))$ and the WGN process $v(k)$.

\subsection{Sparse Gaussian process regression}
\label{subsec:sparse-gp}
A major computational bottleneck of full GPs lies in the construction and inversion of the Gram matrix, which scales as $\mathcal{O}(N^2)$ and $\mathcal{O}(N^3)$, respectively. Additionally, evaluating~\eqref{eq:pred_mean}, and~\eqref{eq:pred_var} requires computing a kernel slice for every data point, resulting in costs of $\mathcal{O}(N)$, and $\mathcal{O}(N^2)$, respectively. Sparse approximations~\cite{JoCa05_UnifiedFramework_Sparse, Ti09_VFE} alleviate this by introducing a smaller set of inducing points yielding a reduced synthetic dataset, $\mathcal{\breve{D}}_{M,i} = \{(\breve{w}_\tau, \breve{z}_{i,\tau} )\}_{\tau=1}^M \approx \mathcal{D}_{N,i}$, with $M \ll N$. This reduces the training complexity to $\mathcal{O}(NM^2)$ and the prediction costs for the mean and variance to $\mathcal{O}(M)$, and $\mathcal{O}(M^2)$, respectively, while ensuring that the resulting predictive distribution $\mathcal{N}\bigl(\breve{\mu}_{z}(w^\ast), \breve{\Sigma}_{z}(w^\ast)\bigr) \approx \mathcal{N}\bigl(\mu_{z}(w^\ast), \Sigma_{z}(w^\ast)\bigr)$.

In this paper, we consider the~\emph{variational free energy} (VFE) method~\cite{Ti09_VFE}, where the goal is to minimize the~\emph{Kullback-Leibler} (KL) divergence between the true GP posterior distribution and an approximated one characterized by $M$ inducing points.

The resulting predictive distribution, i.e., $p(z^\ast_i|\breve{\mathcal{D}}_{M,i}, w^\ast) = \mathcal{N}\bigl(\breve{\mu}_{z,i}(w^\ast), \breve{\sigma}_{z,i}^2(w^\ast)\bigr)$ for each output dimension $z_i^\ast$ with $i \in \mathbb{I}_{1}^{n_\mathrm{z}}$, is described by new predictive mean $\breve{\mu}_{z,i}(w^\ast)$, and variance $\breve{\sigma}_{z,i}^2(w^\ast)$ functions expressed\footnote{Similar to the full GP notation, $\left[ K_{\breve{w}\breve{w},i} \right]_{\tau, \tilde{\tau}} = \kappa_i(\breve{w}_\tau, \breve{w}_{\tilde{\tau}})$, $\left[ K_{\breve{w}w,i} \right]_{\tau, \bar{\tau}} = \kappa_i(\breve{w}_\tau, w_{\bar{\tau}})$, and $\left[ K_{\breve{w},i}(w^\ast) \right]_{\tau} = \kappa_i(\breve{w}_\tau, w^\ast)$, with $\tau, \tilde{\tau} \in \mathbb{I}_{1}^{M}$, and $\bar{\tau} \in \mathbb{I}_{1}^{N}$.} by
\begin{subequations}
\label{eq:sparse_predictive_distr}
\begin{align}
&\breve{\mu}_{z,i}(w^\ast) = K_{\breve{w},i}^\intercal(w^\ast)
\overbrace{\sigma_{v,i}^{-2}\mathcal{S}_{\breve{w}w,i}^{-1}K_{\breve{w}w,i}Z_i}^{\breve{\alpha}_i},
\label{eq:sparse_pred_mean} \\
& \breve{\sigma}_{z,i}^2(w^\ast) = \kappa_i(w^\ast,w^\ast) \notag \\
&\qquad \qquad \; - K_{\breve{w},i}^\intercal(w^\ast)
\bigl(K_{\breve{w}\breve{w},i}^{-1} - \mathcal{S}_{\breve{w}w,i}^{-1}\bigr)
K_{\breve{w},i}(w^\ast), \label{eq:sparse_pred_var}
\end{align}
\end{subequations}
with $\mathcal{S}_{\breve{w}w,i} = K_{\breve{w}\breve{w},i} + \sigma_{v,i}^{-2} K_{\breve{w}w,i} K_{\breve{w}w,i}^\intercal$.

To learn the inducing inputs $\{\breve{w}_\tau\}_{\tau=1}^M$ and hyperparameters $\theta_i$, we minimize the VFE, which is equivalent to maximizing the \emph{evidence lower bound} (ELBO) (see~\cite{Ti09_VFE} for details). The resulting objective is
\begin{multline}
-\log p(Z_i \mid W, \breve{W}, \theta_i)
= \frac{1}{2}\Bigl(
\log \det\bigl(\mathcal{Q}(\breve{\theta}_i) + \gamma^{-1}(\breve{\theta}_i) I_N\bigr) \\
+ Z_i^\intercal \bigl(\mathcal{Q}(\breve{\theta}_i) + \gamma^{-1}(\breve{\theta}_i) I_N\bigr) Z_i
+ \gamma(\breve{\theta}_i)\operatorname{Tr}\bigl(K_{ww,i} - \mathcal{Q}(\breve{\theta}_i)\bigr)
\Bigr),
\label{eq:inducing-opt}
\end{multline}
with $\breve{W} = \mathrm{col}\left(\breve{w}_1,\dots,\breve{w}_M\right)$, $\breve{\theta}_i = \mathrm{col}(\theta_i, \breve{W})$, $\gamma(\breve{\theta}_i) = \sigma_{v,i}^{-2}$, and $\mathcal{Q}(\breve{\theta}_i) = K^\intercal_{\breve{w},w,i}(\breve{\theta}_i) K_{\breve{w}\breve{w},i}^{-1}(\breve{\theta}_i) K_{\breve{w}w,i}(\breve{\theta}_i)$. Then, the optimal hyperparameters $\breve{\theta}_i$ can be chosen as
\begin{equation}
    \breve{\theta}_i^\star = \argmin_{\breve{\theta}_i} \bigl(-\log p(Z_i |W, \breve{W}, \theta_i)\bigr).
\end{equation}

\section{GP-MPC formulation}
\label{sec:mpc}
This section introduces the GP-MPC approach based on the GP-augmented dynamics model~\eqref{eq:aug_model}. We first derive the full stochastic optimal control problem in~\Cref{subsec:stochastic-mpc}, based on the considered control objective. Since solving this problem is generally intractable, see~\cite{HeLiZe18_GPMPC2} and~\Cref{subsec:stochastic-mpc}, we adopt two standard approximation schemes in~\Cref{subsec:model-propagation}:
\begin{enumerate*}[label=(\roman*)]
\item first-order Taylor expansion~\cite{GiRaMu02_Uncertain_GP_Taylor}, and
\item the MM method~\cite{QuGiRa03_Uncertain_GP_EMM}.
\end{enumerate*}
As our first contribution~\ref{contributions:item-one}, we derive closed-form MM predictors for sparse GPs, enabling efficient propagation of state distributions; an extension not previously addressed in the literature. As shown in~\Cref{subsec:nmpc}, both approximations yield tractable deterministic nonlinear programs and scale to large datasets, with MM additionally preserving predictive accuracy.

\subsection{Stochastic predictive controller}
\label{subsec:stochastic-mpc}
As introduced in~\Cref{subsec:sys-model}, we consider a predictive reference tracking problem, now reformulated under the GP-augmented dynamics model~\eqref{eq:aug_model}. Due to the GP-based process model, the predicted state distributions at each discrete time-step $k$ necessitate a probabilistic formulation of the MPC problem. A standard quadratic cost is used, and the objective is to minimize its expected value:
\begin{equation}
    \begin{split}
        J(k) &=  \mathbb{E}\Biggl\{\sum_{i=0}^{N_{\mathrm{p}}} \lVert x(i|k) - r(i|k) \rVert_Q^2 + \sum_{i=0}^{N_{\mathrm{p}}-1} \lVert u(i|k) \rVert_R^2 \Biggr\},
    \end{split}
    \label{eq:stoch_mpc:cost_func}
\end{equation}
where $N_{\mathrm{p}}$ is the prediction horizon, $r(i|k) \in \mathbb{R}^{n_\mathrm{x}}$ denotes the reference trajectory over $N_{\mathrm{p}}$, and $Q \in \mathrm{S}^{n_\mathrm{x}}_{\succ}$, $R \in \mathrm{S}^{n_\mathrm{u}}_{\succ}$ are weighting matrices that shape the trade-off between tracking performance and control effort. In cases where penalizing the deviation of $u(i|k)$ from zero is not compatible with the desired state reference $r(i|k)$, a cost penalizing the control increments $\Delta u(i|k)$ may be used instead, with $\Delta u(i|k)=u(i|k)-u(i-1|k)$ and $u(-1|k)\coloneqq u(k-1)$. For simplicity of exposition, we adopt the cost in~\eqref{eq:stoch_mpc:cost_func} in what follows.

Due to the stochastic nature of the prediction model, probabilistic (chance) constraints are imposed on the states to ensure constraint satisfaction with a given probability, while input constraints remain deterministic. Specifically,
\begin{subequations}\label{eq:chance_constr}
\begin{align}
    &\mathrm{Pr}\left(x(i|k)\in \mathcal{X}\right) \geq p_x, \; \forall i \in  \mathbb{I}_{1}^{N_{\mathrm{p}}} \label{eq:stoch_mpc:state_constr}, \\
    &u(i|k)\in \mathcal{U}, \; \forall i \in  \mathbb{I}_{0}^{N_{\mathrm{p}}-1},\label{eq:stoch_mpc:input_constr}
\end{align}
\end{subequations}
where $p_x \in (0, 1)$ denotes the required probability level for constraint satisfaction.
The sets $\mathcal{X} \subseteq \mathbb{R}^{n_\mathrm{x}}$ and $\mathcal{U} \subseteq \mathbb{R}^{n_\mathrm{u}}$ represent the polytopic state $x$ and input $u$ constraint sets defined as the intersection of $n_{\mathrm{CX}}$ and $n_{\mathrm{CU}}$ number of half-spaces, respectively.
That is
\begin{subequations}
    \begin{align}
        \mathcal{X} &= \bigcap_{j=1}^{n_{\mathrm{CX}}} \Bigl\{x \in \mathbb{R}^{n_\mathrm{x}} \mid \alpha_{x,j}^\intercal x \leq b_{x,j} \Bigr\}, \label{eq:polytopic_x}
        \\
        \mathcal{U} &= \bigcap_{j=1}^{n_{\mathrm{CU}}} \Bigl\{u \in \mathbb{R}^{n_\mathrm{u}} \mid \alpha_{u,j}^\intercal u \leq b_{u,j} \Bigr\}. \label{eq:polytopic_u}
    \end{align}
    \label{eq:init_polytops}
\end{subequations}
Then, the stochastic MPC problem is formulated as
\begin{subequations}
\begin{align}
\min_{\{u(i|k)\}_{i=0}^{N_{\mathrm{p}}-1}} &J(k) = \mathbb{E}\Biggl\{\sum_{i=0}^{N_{\mathrm{p}}} \lVert x(i|k) - r(i|k) \rVert_Q^2  \notag \\
& \qquad \quad \; + \sum_{i=0}^{N_{\mathrm{p}}-1} \lVert u(i|k) \rVert_R^2 \Biggr\},
 \label{eq:stoch_mpc:cost} \\
\mathrm{s.t.} \quad & x(i+1|k) = f\bigl(x(i|k),u(i|k)\bigr) + \hat z(i|k),
    \label{eq:stoch_mpc:pred_model} \\
        &p\bigl(\hat z(i|k)\bigr) = \int \mathcal{N}\Bigl(\hat z(i|k) \mid \breve{\mu}_z\bigl(w(i|k)\bigr), \notag \\
        &\qquad \qquad \breve{\Sigma}_z\bigl(w(i|k)\bigl) \Bigr) p\bigl(w(i|k)\bigr) \mathrm{d}w, 
     \label{eq:stoch_mpc:z_part} \\
    &\mathrm{Pr}\left(x(i+1|k)\in \mathcal{X}\right) \geq p_x, \\
    & u(i|k)\in \mathcal{U}, \\
    & x(0|k) = x(k), \label{eq:stoch_mpc:x0}
\end{align}
\label{eq:stoch_mpc}
\end{subequations}
with $w(i|k) = \mathrm{col}\bigl(x(i|k), u(i|k)\bigr)$.

Here,~\eqref{eq:stoch_mpc:z_part} represents the marginal predictive density of the GP correction term $\hat z(i|k)$ under uncertain input $w(i|k)$, obtained by marginalizing the conditional Gaussian GP posterior with respect to the distribution of $w(i|k)$. Together with the nonlinear nominal term $f\bigl(x(i|k),u(i|k)\bigr)$ in~\eqref{eq:stoch_mpc:pred_model}, it induces the predictive density of the next state $x(i+1|k)$.

As the predictive model~\eqref{eq:stoch_mpc:pred_model} is recursively propagated, it generates a sequence of generally non-Gaussian state distributions. For $i=0$, the argument $w(0|k)$ is deterministic, so that~\eqref{eq:stoch_mpc:z_part} reduces to the Gaussian GP posterior of $\hat z(0|k)$, and thus $x(1|k)$ is Gaussian distributed. However, for $i\geq 1$, the input $w(i|k)$ becomes random due to the uncertainty in past states. Consequently, both the nonlinear mapping induced by $f\bigl(x(i|k),u(i|k)\bigr)$ and the marginalization of $\hat z(i|k)$ in~\eqref{eq:stoch_mpc:z_part} lead, in general, to non-Gaussian predictive state distributions.

As a result, the recursive nature of GP-MPC leads to increasingly complex non-Gaussian distributions, rendering exact stochastic optimization intractable. To address this, in~\Cref{subsec:model-propagation}, we adopt two standard methods for approximating~\eqref{eq:stoch_mpc:pred_model} and~\eqref{eq:stoch_mpc:z_part}:
\begin{enumerate*}[label=(\roman*)]
\item first-order Taylor expansion, and
\item the MM approach.
\end{enumerate*}
For the latter, we further present an extension to sparse GPs according to~\Cref{subsec:sparse-gp}. 

\subsection{Augmented model approximation}
\label{subsec:model-propagation}
Following existing GP-MPC works~\cite{HeKaZe20_GPMPC1, PoPeTo23_LPVGPMPC}, we approximate the predicted state distribution at each step $i \in \mathbb{I}_{1}^{N_{\mathrm{p}}}$ as Gaussian, i.e., $x(i|k) \sim \mathcal{N}\bigl(\mu_x(i|k), \Sigma_x(i|k)\bigr)$. Since the predicted input $u(i|k)$ is deterministic, this induces a Gaussian distribution for $w(i|k)=\mathrm{col}\bigl(x(i|k),u(i|k)\bigr)$. For the sake of readability, we drop the prediction and time indices in the remainder of this subsection. Thus, given $w \sim \mathcal{N}\bigl(\mu_w, \Sigma_w\bigr)$ with $\mu_w = \mathrm{col}(\mu_x, u)$, and $\Sigma_w = \mathrm{diag}(\Sigma_x, 0 I_{n_\mathrm{u}})$, and assuming conditional independence of the nominal model $f$ and the GP correction term $\hat z$ in~\eqref{eq:stoch_mpc:pred_model}, we define the approximate joint distribution of the two random vectors as
\begin{equation}
    \begin{split}
        \begin{bmatrix}
        f(w) \\ \hat{z}
        \end{bmatrix}
        \sim \mathcal{N}\left(
        \begin{bmatrix}
            \mu_f \\ \bar{\mu}_z
        \end{bmatrix},
        \begin{bmatrix}
            \Sigma_{f} & \Sigma_{f,z} \\
            \Sigma_{f,z}^\intercal & \bar\Sigma_z
        \end{bmatrix}\right).
    \end{split}
    \label{eq:state_prop:joint_distr}
\end{equation}
As $f$ is a deterministic nonlinear function, no closed-form expression is generally available for marginalizing its output over a Gaussian input $w$.
In contrast, the GP term $\hat{z}$ permits analytical moment computation due to the Gaussian structure of the posterior. Hence, the moments of $f$, namely $\mu_f$ and $\Sigma_f$, as well as the $f$-dependent part of the cross-covariance $\Sigma_{f,z}$, are approximated through Taylor linearization, whereas the GP term $\hat z$ can be handled either by Taylor approximation or by MM-based propagation.

Next, we detail these two propagation techniques and then specialize them to sparse GP models, including our closed-form MM expressions for the sparse case.
\subsubsection{Propagation with Taylor series}
\label{subsubsec:taylor_series}
Let $\mu_z(w)$ and $\Sigma_z(w)$ denote the predictive mean and covariance of the GP correction term $\hat z$, respectively. Linearizing $f(w)$, $\mu_z(w)$, and $\Sigma_z(w)$ around the mean $\mu_w$ of the Gaussian input $w$ yields the following Gaussian approximations (see~\cite{GiRaMu02_Uncertain_GP_Taylor} for details):
\begin{equation}
    f(w) \sim \mathcal{N}\bigl(\underbrace{f(\mu_w)}_{\mu_f}, \underbrace{\nabla_{w}f(\mu_w) \Sigma_w \nabla^\intercal_{w}f(\mu_w)}_{\Sigma_f} \bigr),
    \label{eq:f_taylor_gaussian}
\end{equation}
\begin{equation}
    \hat{z} \sim \mathcal{N}\bigl(\underbrace{{\mu}_z(\mu_w)}_{\bar{\mu}_z}, \underbrace{{\Sigma}_z(\mu_w) + \nabla_{w}{\mu}_z(\mu_w) \Sigma_w \nabla^\intercal_{w} {\mu}_z(\mu_w)}_{\bar{\Sigma}_z} \bigr).
    \label{eq:gp_taylor_gaussian}
\end{equation}
Furthermore, the cross-covariance between $f(w)$ and $\hat z$ is approximated as
\begin{subequations}
    \label{eq:gp_taylor_crossvar}
    \begin{align}
        \Sigma_{f,z} &= \nabla_{w}f(\mu_w) \bar{\Sigma}_{f,z}, \label{eq:gp_taylor_f_crossvar} \\
        \bar{\Sigma}_{f,z} &= \Sigma_w \nabla^\intercal_{w} {\mu}_z(\mu_w). \label{eq:gp_taylor_gp_crossvar}
    \end{align}
\end{subequations}
The term $\bar{\Sigma}_{f,z}$ in~\eqref{eq:gp_taylor_gp_crossvar} can alternatively be evaluated by MM, as discussed next.

\subsubsection{Propagation with moment matching}
\label{subsubsec:moment_matching}
For GP models with SE kernels $\kappa$ (see~\eqref{eq:se_kernel}), marginalizing the GP predictive density of $\hat z$ in~\eqref{eq:stoch_mpc:z_part} under a Gaussian input $w$ yields closed-form expressions for the GP moments $\bar{\mu}_z$ and $\bar{\Sigma}_z$. In addition, a closed-form expression can be obtained for the GP-dependent term $\bar{\Sigma}_{f,z}$ in~\eqref{eq:gp_taylor_gp_crossvar}, while the contribution of the nominal model $f$ to $\Sigma_{f,z}$ remains Taylor-based through $\nabla_w f(\mu_w)$ in~\eqref{eq:gp_taylor_f_crossvar}.

For full GP models, such MM expressions are well documented in the literature; see, e.g.,~\cite{DeHuHa_AnalyticMM} and~\cite[Section~2.3.2]{De_PILCO_PhD}. Therefore, for the sake of compactness, we do not reproduce these formulas here. However, these results do not directly carry over to sparse GPs, because the corresponding sparse predictive mean and covariance have a different form. This motivates the sparse-GP derivation discussed next.

\subsubsection{Sparse propagation}
\label{subsubsec:sparse_mm}
For sparse GP models, the Taylor approximation in~\eqref{eq:f_taylor_gaussian}--\eqref{eq:gp_taylor_crossvar} retains exactly the same form, with the GP predictive mean and covariance $\mu_z(\cdot)$ and $\Sigma_z(\cdot)$ evaluated using the sparse predictive distribution in~\eqref{eq:sparse_predictive_distr}, i.e., $\breve{\mu}_z(\cdot)$ and $\breve{\Sigma}_z(\cdot)$. By contrast, for MM, the standard full-GP expressions from~\cite{DeHuHa_AnalyticMM} are no longer directly applicable. Therefore, closed-form expressions must be derived for the sparse case in order to obtain $\bar{\mu}_z$, $\bar{\Sigma}_z$, and $\bar{\Sigma}_{f,z}$.

Using the SE kernel $\kappa$ from~\eqref{eq:se_kernel}, these sparse-GP MM expressions are given in~\eqref{eq:gp_emm_gaussian}, with proof provided in the Appendix:
\begin{subequations}\label{eq:gp_emm_gaussian}
    \begin{align}
        & \bar{\mu}_{z,i}(\mu_w, \Sigma_w) = \breve{\alpha}_i^\intercal l_i(\mu_w, \Sigma_w), \label{eq:gp_emm_gaussian_means} \\
        & \bar{\Sigma}_{z,i,i}(\mu_w, \Sigma_w) = \sigma_{i}^2 + \sigma_{v,i}^2 + \breve{\alpha}_i^\intercal L_{i,i}(\mu_w, \Sigma_w) \breve{\alpha}_i \notag \\
            & \qquad \qquad \qquad - \bar{\mu}^2_{z,i}(\mu_w, \Sigma_w) - \mathcal{L}_{\breve{w}w,i}(\mu_w, \Sigma_w),
         \label{eq:gp_emm_gaussian_vars} \\
        & \bar{\Sigma}_{z,i,j}(\mu_w, \Sigma_w) = \breve{\alpha}_i^\intercal L_{i,j}(\mu_w, \Sigma_w) \breve{\alpha}_j \notag \\ 
            & \qquad \qquad \qquad - \bar{\mu}_{z,i}(\mu_w, \Sigma_w) \bar{\mu}_{z,j}(\mu_w, \Sigma_w), \; \forall i\neq j, \label{eq:gp_emm_gaussian_covars} \\
        & \bar{\Sigma}_{f,z,:,i}(\mu_w, \Sigma_w) = \sum_{\tau=1}^{M} \breve{\alpha}_{i,\tau} l_{i,\tau} \Sigma_w (\Lambda_i + \Sigma_w)^{-1} (\breve{w}_\tau - \mu_w), \label{eq:gp_emm_gaussian_crosscovars}
    \end{align}
\end{subequations}
with $\mathcal{L}_{\breve{w}w,i}(\mu_w, \Sigma_w) = \mathrm{Tr}\bigl(( K_{{\breve{w}}{\breve{w}},i}^{-1} - \mathcal{S}_{\breve{w}w,i}^{-1} ) L_{i,i}(\mu_w, \Sigma_w)  \bigr)$, for each $i, j \in \mathbb{I}_{1}^{n_{\mathrm{z}}}, \tau, \tilde\tau \in \mathbb{I}_{1}^{M}$. The auxiliary terms $l_{i,\tau}$ and $\left[L_{i,j} \right]_{\tau, \tilde{\tau}}$ are defined below using the kernel in~\eqref{eq:kernel}:
\begin{subequations}
    \label{eq:gp_emm_auxiliary}
    \begin{align*}
        \begin{split}
            \begin{aligned}[t] & l_{i,\tau}(\mu_w, \Sigma_w) = \sigma_{i}^2 \mathrm{det}\bigl(\Lambda_i^{-1}\Sigma_w + I\bigr)^{-\frac{1}{2}} \\
            & \qquad \exp\{-\frac{1}{2}\left(\mu_w - \breve{w}_\tau \right)^\intercal\left(\Lambda_i + \Sigma_w \right)^{-1}\left(\mu_w - \breve{w}_\tau\right)\}, \end{aligned} \noindent
        \end{split} \\  
        \begin{split}
            \begin{aligned}[t] & \left[L_{i,j} \right]_{\tau, \tilde{\tau}}(\mu_w, \Sigma_w) = \sigma_{i}^2 \sigma_{j}^2 \mathrm{det} \bigl( ( \Lambda_i^{-1} + \Lambda_j^{-1} ) \Sigma_w + I\bigr)^{-\frac{1}{2}} \\
                & \qquad \exp \{-\frac{1}{2}\left(\breve{w}_\tau - \breve{w}_{\tilde\tau}\right)^\intercal \left(\Lambda_i +\Lambda_j\right)^{-1}\left(\breve{w}_\tau - \breve{w}_{\tilde\tau} \right) \}  \\
                & \qquad \exp \{-\frac{1}{2} \left(q_{i,j} - \mu_w \right)^\intercal F^{-1}\left(q_{i,j} - \mu_w \right) \}, \end{aligned}
        \end{split} 
    \end{align*}
\end{subequations}
where $q_{i,j} = \Lambda_j(\Lambda_i + \Lambda_j)^{-1}\breve{w}_\tau + \Lambda_i(\Lambda_i + \Lambda_j)^{-1}\breve{w}_{\tilde\tau}$, and $F = \bigl((\Lambda_i^{-1}+\Lambda_j^{-1})^{-1} + \Sigma_w\bigr)$.
\begin{remark}
The closed-form MM expressions in~\eqref{eq:gp_emm_gaussian} differ from the standard full GP MM formulation~\cite{DeHuHa_AnalyticMM} in several key aspects. First, the trace correction term in the variances $\bar{\Sigma}_{z,i,i}$~\eqref{eq:gp_emm_gaussian_vars} uses $\mathcal{L}_{\breve{w}w,i} = \mathrm{Tr}\bigl(( K_{{\breve{w}}{\breve{w}},i}^{-1} - \mathcal{S}_{\breve{w}w,i}^{-1}) L_{i,i} \bigr)$, where $\mathcal{S}_{\breve{w}w,i}$ is the variational posterior precision, replacing the full GP counterpart $\mathrm{Tr}(K_{ww}^{-1} L_{i,i})$. Second, the predictive mean, covariance, and cross-covariance are computed using sparse GP parameters $\breve{\alpha}_i$ and inducing inputs $\breve{w}$ rather than their full GP equivalents. Lastly, the additive term $\sigma_{v,i}^2$ in the diagonal of the predictive covariance~\eqref{eq:gp_emm_gaussian_vars} arises from the WGN term in the kernel definition~\eqref{eq:kernel}, which explicitly models process noise in the system dynamics.
\end{remark}

\subsubsection{Approximated prediction model}
Utilizing the approximated joint distribution~\eqref{eq:state_prop:joint_distr} of $f(w)$ and $\hat z$, where the quantities $\mu_f$, $\bar{\mu}_z$, $\Sigma_f$, $\bar{\Sigma}_z$, and $\Sigma_{f,z}$ depend on the chosen approximation method (i.e., Taylor or MM), the generally non-Gaussian stochastic prediction model~\eqref{eq:stoch_mpc:pred_model}--\eqref{eq:stoch_mpc:z_part} is approximated by a tractable Gaussian recursion in which only the first two moments are propagated. Accordingly, the predictive state distribution is approximated as $x(i|k) \sim \mathcal{N}\bigl(\mu_x(i|k), \Sigma_x(i|k)\bigr)$, with moments given by
\begin{subequations}
    \label{eq:aug_model_compact_form}
    \begin{align}
        &\mu_{x}(i+1|k) = f\bigl(\mu_{w}(i|k)\bigr) + \bar{\mu}_{z}\bigl(\mu_{w}(i|k), \Sigma_{w}(i|k)\bigr), \label{eq:aug_model_compact_form:mean} \\
        &\Sigma_{x}(i+1|k) = \bar{A}\bigl(\mu_{w}(i|k)\bigr) \notag  \\
& \; \begin{bmatrix}
        \Sigma_{w}(i|k) & \bar\Sigma_{f,z} \bigl(\mu_{w}(i|k), \Sigma_{w}(i|k)\bigr) \\
        \star           & \bar\Sigma_z \bigl( \mu_{w}(i|k), \Sigma_{w}(i|k) \bigr)
    \end{bmatrix} \bar{A}^\intercal \bigl(\mu_{w}(i|k)\bigr),
        \label{eq:aug_model_compact_form:covar}
    \end{align}
\end{subequations}
with $\bar{A}\bigl(\mu_{w}(i|k)\bigr) = \begin{bmatrix} \nabla_{w}f\bigl(\mu_{w} (i|k) \bigr) & I\end{bmatrix}$ and $i \in \mathbb{I}_{0}^{N_{\mathrm{p}}}$.

\subsection{Nonlinear predictive control problem}
\label{subsec:nmpc}
The approximated Gaussian distribution of the state in~\eqref{eq:aug_model_compact_form} enables efficient propagation of the prediction model over the horizon $N_{\mathrm{p}}$, and allows reformulating the stochastic MPC problem~\eqref{eq:stoch_mpc} into a deterministic nonlinear program. This is achieved by computing the expected cost and approximating the state chance constraints in a deterministic manner. The expected cost~\eqref{eq:stoch_mpc:cost_func} becomes
\begin{equation}
\begin{aligned}[b]
     \bar{J}(k) &= \sum_{i=0}^{N_{\mathrm{p}}} \Bigl( \lVert \mu_x(i|k) - r(i|k) \rVert_Q^2 + \mathrm{Tr} \bigl( Q \Sigma_x(i|k) \bigr) \Bigr) \\ &\quad   + \sum_{i=0}^{N_{\mathrm{p}}-1} \lVert u(i|k) \rVert_R^2. \label{eq:cost_func_aprox}
\end{aligned}
\end{equation}
The state chance constraints~\eqref{eq:stoch_mpc:state_constr} are tightened on the state mean $\mu_x(i|k)$ using \emph{probabilistic reachable sets} (PRS) on their respective errors (see~\cite{HeKaZe20_GPMPC1} for details) yielding
\begin{equation}
        \mu_x(i|k) \in \mathcal{Z}\bigl(\Sigma_x(i|k)\bigr), \label{eq:deter_mpc:state_mean_constr}
\end{equation}
where $\mathcal{Z}\bigl(\Sigma_x(i|k)\bigr) \subseteq\mathbb{R}^{n_\mathrm{x}}$ is the polytopic constraint set defined as the intersection of the redefined half-spaces $n_{\mathrm{CX}}$. The redefined half-spaces are reformulated as
\begin{equation}
\label{eq:polytopic_z}
\begin{aligned}[b]
\mathcal{Z}\bigl(\Sigma_x(i|k)\bigr)
= \bigcap_{j=1}^{n_{\mathrm{CX}}} & \Bigl\{ x\in \mathbb{R}^{n_\mathrm{x}}
\mid \alpha_{x,j}^\intercal x \\
& \qquad \leq b_{x,j} - c_{x,j}\bigl(\Sigma_x(i|k)\bigr) \Bigr\},
\end{aligned}
\end{equation}
where $c_{x,j}\bigl(\Sigma_x(i|k)\bigr)= \Phi^{-1}(p_x)\sqrt{\alpha_{x,j}^\intercal\Sigma_x(i|k)\alpha_{x,j}}$, with $\Phi^{-1}(\cdot)$ being the inverse cumulative distribution function of the standard Gaussian distribution.
Finally, the current GP-MPC formulation is based on the assumption that the full noise-free state $x(k) =x(0|k)$ is measurable at every control cycle, which translates to
\begin{equation}
    \mu_{x}(0|k) = x(k), \qquad \Sigma_{x}(0|k) = 0I_{n_\mathrm{x}}.
    \label{eq:nmpc:init_cond}
\end{equation}

Based on these, the stochastic NMPC~\eqref{eq:stoch_mpc} can be formulated as a deterministic NMPC in terms of
\begin{subequations}
\begin{align}
        & \min_{\{u(i|k)\}_{i=0}^{N_{\mathrm{p}}-1}} \bar{J}(k) = \sum_{i=0}^{N_{\mathrm{p}}} \Bigl( \lVert \mu_x(i|k) - r(i|k) \rVert_Q^2 \notag \\
& \qquad \qquad + \mathrm{Tr} \bigl( Q \Sigma_x(i|k) \bigr) \Bigr) + \sum_{i=0}^{N_{\mathrm{p}}-1} \lVert u(i|k) \rVert_R^2,
    \label{eq:nmpc:cost} \\
&\mathrm{s.t.} \nonumber \\
&\mu_{x}(i+1|k) = f\bigl(\mu_{w}(i|k)\bigr) + \bar{\mu}_{z}\bigl(\mu_{w}(i|k), \Sigma_{w}(i|k)\bigr),
\label{eq:nmpc:pred_model_mean} \\
&\Sigma_{x}(i+1|k) = \bar{A}\bigl(\mu_{w}(i|k)\bigr) \notag  \\
& \; \begin{bmatrix}
        \Sigma_{w}(i|k) & \bar\Sigma_{f,z} \bigl(\mu_{w}(i|k), \Sigma_{w}(i|k)\bigr) \\
        \star           & \bar\Sigma_z \bigl( \mu_{w}(i|k), \Sigma_{w}(i|k) \bigr)
    \end{bmatrix} \bar{A}^\intercal \bigl(\mu_{w}(i|k)\bigr),
\label{eq:nmpc:pred_model_var} \\
& \mu_x(i|k) \in \mathcal{Z}\bigl(\Sigma_x(i|k)\bigr)
\\
& u(i|k)\in \mathcal{U},   
\\
& \mu_{x}(0|k) = x(k), \quad \Sigma_{x}(0|k) = 0I_{n_\mathrm{x}}. \label{eq:nmpc:x0}
\end{align}
\label{eq:nmpc}
\end{subequations}
The resulting nonlinear program can be solved to local optima by interior-point~\cite{INTERIOR_POINT} or \emph{sequential QP} (SQP)~\cite{SQP1} methods.

\section{Surrogate LPV form-based solution}
\label{sec:lpv-sol}
While the nonlinear optimization problem in~\eqref{eq:nmpc} yields a locally optimal solution to the GP-MPC problem, the strong nonlinearities in~\eqref{eq:aug_model_compact_form} render the optimization computationally expensive.~\cite{PoPeTo23_LPVGPMPC} introduces an accelerated solution based on LPV iterations, where a series of QPs is solved until convergence. However, the LPV formulation in~\cite{PoPeTo23_LPVGPMPC} does not retain the exact nonlinear model and the variance dynamics remain fixed with respect to the introduced scheduling variables. In this work, we adopt an LPV conversion technique based on the~\emph{fundamental theorem of calculus} (FTC)~\cite{LPV-FTC-new}, which yields an exact, approximation-free affine representation of~\eqref{eq:aug_model_compact_form}. This enables transforming the nonlinear optimization problem into an iterative sequence of QPs while maintaining feasible and expressive variance propagation
\footnote{The proposed LPV approach relies on an FTC embedding that preserves the nonlinear dynamics in affine form, whereas SQP is derived from local linearization of the nonlinear optimization problem to update the Newton search direction. Recent work~\cite{SQP-LPVMPC} shows that, under specific anchor-point choices, an FTC LPV-MPC step can coincide with the SQP subproblem. While SQP benefits from well-established local convergence guarantees under mild assumptions, such guarantees have not been established in the present work for the proposed iterative LPV scheme, which instead prioritizes computational tractability and fidelity to the underlying augmented nonlinear model.}.

\subsection{LPV predictive control formulation}
To provide an efficient auto-conversion method of the GP-based augmented state-space model~\eqref{eq:aug_model} to an LPV state-space equivalent, we utilize the FTC method, which acts as a function factorization. Given a continuously differentiable function $g: \mathbb{R}^n \rightarrow \mathbb{R}^m$, it holds for all $\varsigma, \tilde{\varsigma} \in \mathbb{R}^n$ that~\cite[Appendix~C]{Ko_Phd_FTCLPV}
\begin{equation}
    g(\varsigma) - g(\tilde{\varsigma}) = \left( \int_0^1 \dfrac{\partial g}{\partial \varsigma} \bigl( \tilde{\varsigma} + \lambda ( \varsigma - \tilde{\varsigma} ) \bigr) \mathrm{d}\lambda \right) (\varsigma - \tilde{\varsigma}).
    \label{eq:ftc_def}
\end{equation}
Here, the Jacobian $\frac{\partial g}{\partial \varsigma}\bigl(\tilde{\varsigma}+\lambda(\varsigma-\tilde{\varsigma})\bigr)\in\mathbb{R}^{m\times n}$ is evaluated along the straight-line segment connecting the anchor point $\tilde{\varsigma}$ to the query point $\varsigma$, and the integral is taken component-wise along this path.

Without loss of generality, let $\vartheta$ and $\zeta$ denote the nonlinear mappings of the predictive mean and covariance from~\eqref{eq:aug_model_compact_form:mean} and~\eqref{eq:aug_model_compact_form:covar}, respectively. That is,
\begin{subequations}
    \label{eq:aug_model_general_form}
    \begin{align}
        \mu_{x}(k+1) = \vartheta \bigl(\mu_{w}(k), \Sigma_{w}(k) \bigr), \\
        \Sigma_{x}(k+1) = \zeta \bigl(\mu_{w}(k), \Sigma_{w}(k) \bigr), \label{eq:aug_model_general_form:var}
    \end{align}
\end{subequations}
with $\mu_w = \mathrm{col}(\mu_x, u)$, and $\Sigma_w = \mathrm{diag}(\Sigma_x, 0 I_{n_\mathrm{u}})$.

Provided that $\vartheta$ and $\zeta$ are differentiable, for each time instant $k$ and prediction step $i \in \mathbb{I}_{0}^{N_{\mathrm{p}}-1}$, define
\[
\varsigma(i|k) = \mathrm{col}\bigl(\mu_{x}(i|k), u(i|k), \operatorname{vec}(\Sigma_{x}(i|k))\bigr),
\]
and choose the anchor point as
\[
\tilde{\varsigma}(k) = \mathrm{col}\bigl(x(k), u(k-1), \operatorname{vec}(0I_{n_{\mathrm{x}}})\bigr),
\]
where $x(k)$ denotes the measured state at the current sampling instant and $u(k-1)$ the previously applied input. Then, the map $\vartheta$ can be factorized using the FTC as
\begin{equation*}
\begin{split}
        &\vartheta(\varsigma(i|k)) = \vartheta(\tilde{\varsigma}(k)) \\
        & + \underbrace{\int_0^1 \dfrac{\partial \vartheta}{\partial \mu_{x}} \bigl( \tilde{\varsigma}(k) + \lambda ( \varsigma(i|k) - \tilde{\varsigma}(k) ) \bigr) \mathrm{d}\lambda}_{A_{\vartheta}(\varsigma(i|k))} (\mu_{x}(i|k) - x(k)) \\ 
        & + \underbrace{\int_0^1 \dfrac{\partial \vartheta}{\partial u} \bigl( \tilde{\varsigma}(k) + \lambda ( \varsigma(i|k) - \tilde{\varsigma}(k) ) \bigr) \mathrm{d}\lambda}_{B_{\vartheta}(\varsigma(i|k))} (u(i|k) - u(k-1)) \\
        & + \underbrace{\int_0^1 \dfrac{\partial \vartheta}{\partial \operatorname{vec}(\Sigma_{x})} \bigl( \tilde{\varsigma}(k) + \lambda ( \varsigma(i|k) - \tilde{\varsigma}(k) ) \bigr) \mathrm{d}\lambda}_{C_{\vartheta}(\varsigma(i|k))} \operatorname{vec}(\Sigma_{x}(i|k)),
\end{split}
\end{equation*}
where $\vartheta(\tilde{\varsigma}(k))$ is an affine term, and $\zeta(\varsigma(i|k))$ is similarly defined.

\begin{remark}
The \emph{anchor} point $\tilde{\varsigma}(k)$ serves as a fixed reference for evaluating the nonlinear system dynamics. When applying the FTC LPV auto-conversion approach, the difference between $\vartheta(\varsigma(i|k))$ and $\vartheta(\tilde{\varsigma}(k))$ is computed along the line segment connecting $\varsigma(i|k)$ and $\tilde{\varsigma}(k)$, where the anchor point $\tilde{\varsigma}(k)$ keeps the factorization close to the current state. Choosing $x(k)$ in $\tilde{\varsigma}(k)$ aligns the integral of the Jacobian closely with the actual dynamics, since it shortens the line over which the Jacobian is integrated. Notably, this LPV formulation exactly reproduces the original nonlinear model; the only source of approximation may stem from numerical integration of the Jacobian instead of an analytic solution.
\end{remark}

Introducing the scheduling variable at prediction step $i \in \mathbb{I}_{0}^{N_{\mathrm{p}}-1}$ as
\begin{equation}
    \rho(i|k) = \mathrm{col}\bigl( \mu_{x}(i|k), u(i|k), \operatorname{vec}(\Sigma_{x}(i|k)) \bigr),
    \label{eq:scheduling_vars}
\end{equation}
and using the identity function as the scheduling map, i.e., $\phi\bigl(\rho(i|k)\bigr) = \rho(i|k)$, the FTC-based LPV prediction model of~\eqref{eq:aug_model_general_form} becomes
\begin{subequations}
    \label{eq:lpv:dynamics}
    \begin{align}
        &\mu_{x}(i+1|k) = A_{\vartheta} \bigl(\rho(i|k)\bigr) \bigl(\mu_x(i|k) - x(k)\bigr) \notag \\
        & \quad \quad + B_{\vartheta} \bigl(\rho(i|k)\bigr) \bigl(u(i|k) - u(k-1)\bigr) \notag \\
        & \quad \quad + C_{\vartheta} \bigl(\rho(i|k)\bigr) \operatorname{vec}(\Sigma_{x}(i|k)) + \vartheta\bigl(\tilde{\varsigma}(k)\bigr),
        \label{eq:lpv:mean}
        \\
        &\operatorname{vec}(\Sigma_{x}(i+1|k)) = A_{\zeta} \bigl(\rho(i|k)\bigr) \bigl(\mu_x(i|k) - x(k)\bigr) \notag \\
        & \qquad + B_{\zeta} \bigl(\rho(i|k)\bigr) \bigl(u(i|k) - u(k-1)\bigr) \notag \\
        & \qquad + C_{\zeta} \bigl(\rho(i|k)\bigr) \operatorname{vec}(\Sigma_{x}(i|k)) + \operatorname{vec}\bigl(\zeta(\tilde{\varsigma}(k))\bigr),
        \label{eq:lpv:var}
    \end{align}
\end{subequations}
where $\vartheta(\tilde{\varsigma}(k))$ and $\zeta(\tilde{\varsigma}(k))$ are the constant terms obtained by evaluating the nonlinear maps~\eqref{eq:aug_model_general_form} at the anchor point $\tilde{\varsigma}(k)$.

Although calculating the Jacobian of $\vartheta, \zeta$ with respect to $\mu_x, u, \operatorname{vec}(\Sigma_x)$ and integrating them to obtain the corresponding matrices can be a computationally intensive task, 
efficient algorithmic differentiation tools, like CasADi~\cite{casadi}, and numerical integration algorithms, such as Simpson's $1/3$ or $3/8$ rule, can be efficiently used for this purpose.

The surrogate LPV formulation of the NMPC~\eqref{eq:nmpc} can be completed by reformulating the nonlinear state constraints defined in~\eqref{eq:deter_mpc:state_mean_constr} into
\begin{equation}
    \mu_x(i|k) \in \mathcal{Z}\bigl(\rho(i|k)\bigr).
    \label{eq:lpv:state_mean_constr}
\end{equation}

Then, the stochastic NMPC~\eqref{eq:stoch_mpc} can be formulated into a LPV-MPC form as
\begin{subequations}
\label{eq:lpv-mpc}
\begin{align}
& \min_{\{u(i|k)\}_{i=0}^{N_{\mathrm{p}}-1}} \bar{J}(k) = \sum_{i=0}^{N_{\mathrm{p}}} \Bigl( \lVert \mu_x(i|k) - r(i|k) \rVert_Q^2 \notag \\
& \qquad \qquad + \mathrm{Tr} \bigl( Q \Sigma_x(i|k) \bigr) \Bigr) + \sum_{i=0}^{N_{\mathrm{p}}-1} \lVert u(i|k) \rVert_R^2,
\label{eq:lpv-mpc:cost}\\
& \mathrm{s.t.} \notag \\
& \mu_{x}(i+1|k) = A_{\vartheta}\bigl(\rho(i|k)\bigr)\bigl(\mu_x(i|k)-x(k)\bigr) \notag\\
&\quad + B_{\vartheta}\bigl(\rho(i|k)\bigr)\bigl(u(i|k)-u(k-1)\bigr) \notag\\
&\quad + C_{\vartheta}\bigl(\rho(i|k)\bigr)\operatorname{vec}(\Sigma_x(i|k)) \notag\\
&\quad + \vartheta\bigl(x(k),u(k-1),0I_{n_\mathrm{x}}\bigr),
\label{eq:lpv-mpc:pred_model_mean}\\
& \operatorname{vec}(\Sigma_x(i+1|k)) = A_{\zeta}\bigl(\rho(i|k)\bigr)\bigl(\mu_x(i|k)-x(k)\bigr) \notag\\
&\quad + B_{\zeta}\bigl(\rho(i|k)\bigr)\bigl(u(i|k)-u(k-1)\bigr) \notag\\
&\quad + C_{\zeta}\bigl(\rho(i|k)\bigr)\operatorname{vec}(\Sigma_x(i|k)) \notag\\
&\quad + \operatorname{vec}\bigl(\zeta\bigl(x(k),u(k-1),0I_{n_\mathrm{x}}\bigr)\bigr),
\label{eq:lpv-mpc:pred_model_var}\\
& \mu_x(i|k) \in \mathcal{Z}\bigl(\rho(i|k)\bigr), \label{eq:lpv-mpc:state_mean_constr} \\
& u(i|k)\in \mathcal{U},\\
& \mu_x(0|k)=x(k),\quad \Sigma_x(0|k)=0I_{n_\mathrm{x}}.
\label{eq:lpv-mpc:x0}
\end{align}
\end{subequations}

\begin{remark}
At prediction step $i = 0$, the term $\mu_x(0|k) - x(k)$ in~\eqref{eq:lpv-mpc:pred_model_mean},~\eqref{eq:lpv-mpc:pred_model_var} vanishes, which is consistent with anchoring the FTC-based factorization at the most recent state, preserving the exact nonlinear model up to numerical integration error. The dynamics at this step are still preserved through the affine term $\vartheta(\tilde{\varsigma}(k))$, which captures the nonlinear model evaluated at the anchor point.
\end{remark}

\subsection{Iterative GP-MPC method}
\label{subsec:lpv-gpmpc}
The optimization problem in~\eqref{eq:lpv-mpc} reduces the nonlinear GP-MPC problem in~\eqref{eq:nmpc} to a parametric QP, solvable by standard QP solvers for a given scheduling sequence $\{\rho(i|k)\}_{i=0}^{N_{\mathrm{p}}-1}$ at any time instant $k \in \mathbb{Z}$. If the optimal trajectories $\{\mu_x^{\star}(i|k), \operatorname{vec}(\Sigma_x^{\star}(i|k))\}_{i=1}^{N_{\mathrm{p}}}$, and $\{u^{\star}(i|k)\}_{i=0}^{N_{\mathrm{p}}-1}$ of the nonlinear problem~\eqref{eq:nmpc} are known, they induce an associated scheduling sequence $\{\rho^{\star}(i|k)\}_{i=0}^{N_{\mathrm{p}}-1}$ for which the LPV model~\eqref{eq:lpv:dynamics} exactly matches the nonlinear moment propagation model~\eqref{eq:aug_model_compact_form} along that trajectory. Consequently, the solution of~\eqref{eq:nmpc} is also a valid solution of the LPV-MPC problem in~\eqref{eq:lpv-mpc}.

In general, however, the optimal scheduling sequence $\{\rho^{\star}(i|k)\}_{i=0}^{N_{\mathrm{p}}-1}$ is not available a priori, and~\eqref{eq:lpv-mpc} only approximates~\eqref{eq:nmpc} for other scheduling trajectories. To mitigate this, an iterative procedure is adopted, wherein~\eqref{eq:lpv-mpc} is repeatedly solved by updating the scheduling trajectory based on the resulting solution. This procedure is summarized in~\cref{alg:iterated-lpv-mpc}.

At each time step $k \in \mathbb{Z}$,~\cref{alg:iterated-lpv-mpc} takes as input the state measurement $x(k)$ and the previously converged input sequence $\{u(i|k-1)\}_{i=0}^{N_{\mathrm{p}}-1}$. This sequence is shifted forward by one step to initialize the current sequence $\{u(i|k)\}_{i=0}^{N_{\mathrm{p}}-1}$ as
\begin{equation}
    \label{eq:input_shift}
    u(i|k) =
        \begin{cases}
            u(i+1|k{-}1), & \text{for } i \in \mathbb{I}_0^{N_{\mathrm{p}}-2}, \\
            u(N_{\mathrm{p}}{-}1|k{-}1), & \text{for } i = N_{\mathrm{p}}{-}1.
        \end{cases}
\end{equation}

Then, the moment propagation model~\eqref{eq:aug_model_general_form} is simulated using the state measurement $x(k)$ and the updated input sequence $\{u(i|k)\}_{i=0}^{N_{\mathrm{p}}-1}$ to obtain the predicted state mean and covariance trajectories $\{\mu_x(i|k), \operatorname{vec}(\Sigma_x(i|k))\}_{i=1}^{N_{\mathrm{p}}-1}$. These are used to initialize the scheduling variables via~\eqref{eq:scheduling_vars}.
At initialization (i.e., $k=0$), we define $\{u(i|-1)= 0\}_{i=0}^{N_{\mathrm{p}}-1}$, and so $\rho(i|k)$ is initialized with $\{\mu_x(i|0) = x(0), u(i|0)= 0, \operatorname{vec}(\Sigma_x(i|0)) = 0\}_{i=0}^{N_{\mathrm{p}}-1}$ via~\eqref{eq:scheduling_vars}.

Given the scheduling trajectory, the matrices in~\eqref{eq:lpv:dynamics}, and the polytopic set $\mathcal{Z}$ in~\eqref{eq:lpv:state_mean_constr} are parameterized enabling the solution $\{u(i|k)\}_{i=0}^{N_{\mathrm{p}}-1}$ of the LPV-MPC problem~\eqref{eq:lpv-mpc} by a standard QP. This is used to simulate~\eqref{eq:aug_model_general_form} in order to compute the state variation $\{\mu_x(i|k), \operatorname{vec}(\Sigma_x(i|k))\}_{i=1}^{N_{\mathrm{p}}}$. Then, the scheduling sequence $\rho(i|k)$ is updated as in~\eqref{eq:scheduling_vars} for all $i \in \mathbb{I}_{0}^{N_{\mathrm{p}}-1}$.

This procedure is repeated until the computed scheduling sequence has converged or maximum iterations have been reached. The convergence criterion proposed here is based on the $\ell_\infty$-norm, defined as
\begin{equation}
    \lVert \rho^{(j+1)}(i|k) - \rho^{(j)}(i|k) \rVert_{\infty} \leq \epsilon_{\mathrm{lpv}},
    \label{eq:convergence_criterion}
\end{equation}
where $\epsilon_\mathrm{lpv}$ is a threshold value, and $\rho^{(j)}(i|k)$ is the scheduling sequence at iteration $j$ of the inner loop at time instant $k$.

Upon convergence, $u(0|k)$ is applied to~\eqref{eq:aug_model_general_form} and the procedure is repeated for the next time instance $k+1$. 
\begin{algorithm}
\caption{Iterated LPV-MPC solution of~\eqref{eq:nmpc}}\label{alg:iterated-lpv-mpc}
\begin{algorithmic}[1]
\Require $x(k), \{u(i|k-1)\}_{i=0}^{N_{\mathrm{p}}-1}$
\Ensure $u(0|k)$
\Initialize{set $j \leftarrow 0$, $\mu_x(0|k) \leftarrow x(k)$, simulate~\eqref{eq:aug_model_general_form} \\[]
\quad with $\mu_x(0|k)$ and~\eqref{eq:input_shift} to set $\{\rho^{(j)}(i|k)\}_{i=0}^{N_{\mathrm{p}}{-}1}$ via~\eqref{eq:scheduling_vars} \label{alg:init}} 
\Repeat
    \State solve~\eqref{eq:lpv-mpc} to obtain $u(i|k), \forall i \in \mathbb{I}_{0}^{N_{\mathrm{p}}-1}$ \label{alg:solve}
    \State simulate~\eqref{eq:aug_model_general_form} with $\mu_x(0|k), u(i|k), \forall i \in \mathbb{I}_{0}^{N_{\mathrm{p}}-1}$ \\[]  
    \qquad \quad to obtain $\{\mu_x(i|k), \operatorname{vec}(\Sigma_x(i|k))\}_{i=1}^{N_{\mathrm{p}}}$ \label{alg:simulate}
    \State update $\{\rho^{(j+1)}(i|k)\}_{i=0}^{N_{\mathrm{p}}-1}$ as in~\eqref{eq:scheduling_vars} \label{alg:update}
    \State set $j \leftarrow j + 1$ \label{alg:j_update}
\Until{\eqref{eq:convergence_criterion} or max iterations reached}
\end{algorithmic}
\end{algorithm}

\begin{remark}
    A~\emph{real-time iteration} (RTI) variant of~\cref{alg:iterated-lpv-mpc} executes only one iteration per time step $k$: solving the QP once (\cref{alg:solve}) using the initial scheduling trajectory (\cref{alg:init}), applying $u(0|k)$, and skipping the simulation (\cref{alg:simulate}) and update (\cref{alg:update}) steps. At the next time step $k+1$, a new scheduling trajectory is initialized via~\cref{alg:init}, ensuring real-time feasibility with reduced computational cost.
\end{remark}

\subsection{Precomputed covariance propagation}
\label{subsec:precomputed-cov}
Although full covariance propagation enhances the fidelity of GP-MPC by explicitly accounting for model uncertainty, it imposes a computational burden, particularly for high-dimensional systems. In such cases, the iterative LPV-MPC solution can be reformulated by adopting a precomputed covariance propagation strategy. Rather than including the covariance sequence as decision variables in the QP, the propagation is fixed based on the solution from the previous iteration. This reformulation significantly reduces the number of optimization variables, while still allowing state chance constraints and mean dynamics to be evaluated using the fixed covariance sequence. Following~\cite{HeLiZe18_GPMPC2}, this approach is adapted to the LPV-MPC setting and integrates directly into the QP formulation without additional computational cost.

We incorporate this strategy into the LPV-MPC formulation~\eqref{eq:lpv-mpc} by dropping the LPV-based covariance propagation~\eqref{eq:lpv-mpc:pred_model_var}, which reduces the QP size by removing covariance-related decision variables. Since the covariance sequence is already simulated in the iterative LPV-MPC procedure, specifically during the initialization (\cref{alg:init}) and update (\cref{alg:update}) steps, the precomputed sequence $\{\operatorname{vec}(\Sigma_x(i|k))\}_{i=0}^{N_{\mathrm{p}}-1}$ becomes part of the scheduling variables $\{\rho(i|k)\}_{i=0}^{N_{\mathrm{p}}-1}$. Consequently, the state chance constraints~\eqref{eq:lpv-mpc:state_mean_constr} and mean dynamics~\eqref{eq:lpv-mpc:pred_model_mean} remain satisfied without additional modifications to~\eqref{eq:lpv-mpc}.

This approximation introduces no additional computational overhead, making it well-suited for applications requiring fast control updates. While it leads to a more conservative policy, as in~\cite{HeLiZe18_GPMPC2}, it is advantageous in regions of low model uncertainty or when responsiveness is prioritized over precision.

\section{Application to the Crazyflie 2.1}
\label{sec:application}

In this section we demonstrate on the Crazyflie~2.1 the performance and computational advantages of the accelerated LPV-MPC introduced in~\cref{subsec:lpv-gpmpc} over the standard NMPC solution of the GP-MPC problem in~\cref{subsec:nmpc}. We validate this claim in both high-fidelity simulation and real-world experiments. The simulation model mirrors the Crazyflie~2.1 dynamics and limits and is used to test the controllers under reference trajectories corresponding to aggressive maneuvers and external disturbances. The real-world experiments corroborate real-time feasibility and robustness in an off-board control setup. Compared to prior data-driven MPC for quadrotors~\cite{ToKaFo21_GPMPCDrone}, which treats the GP-augmented model deterministically, we propagate GP uncertainty and couple it into the mean dynamics, constraints, and a variance-dependent objective.

\subsection{Quadrotor modeling}
\label{subsec:quadrotor}
First, we summarize the rigid-body dynamics and the thrust/torque allocation used throughout.
Let $\mathcal{F}^i$ be the inertial frame in \emph{East-North-Up} (ENU) convention, and $\mathcal{F}^b$ be the body-fixed frame. With a ZYX (yaw-pitch-roll) Euler parametrization $(\psi,\theta,\phi)$ to characterize three consecutive rotations in $z, y$, and $x$ body axes, the rotation from $\mathcal{F}^b$ to $\mathcal{F}^i$ is $R^i_b(\phi,\theta,\psi)\in\mathrm{SO}(3)$, with $\mathrm{SO}(3)$ being the special orthogonal group.
\begin{figure}[t]
  \centering
   \includegraphics[width=0.6\linewidth]{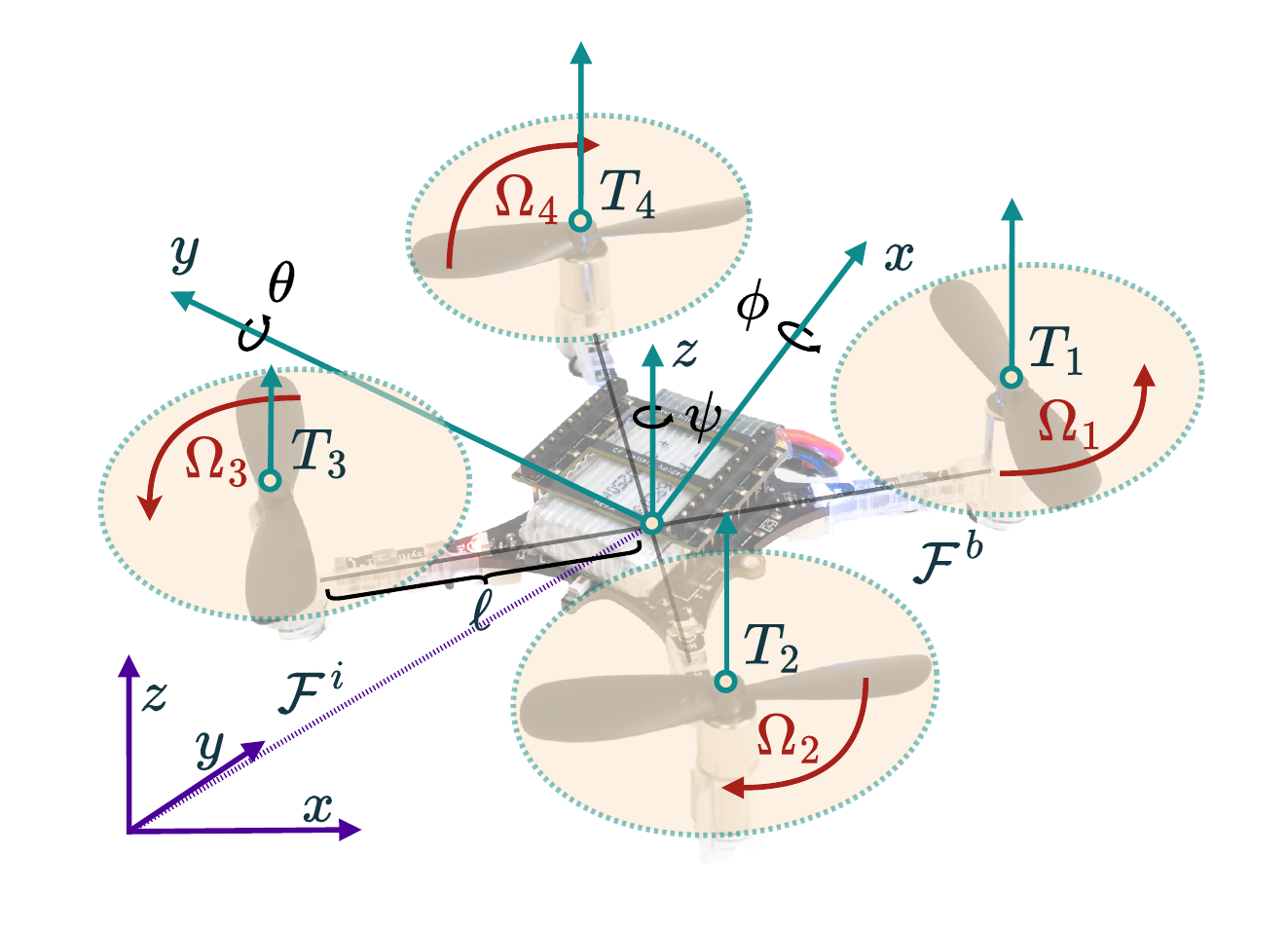}
   \caption{Perspective view of the geometric relationship between $\mathcal{F}^i$, and $\mathcal{F}^b$, as well as the direction of the rotor angular velocities $\Omega_j$ and thrusts $T_j$ for $j \in \mathbb{I}_i^4$.}
   \label{fig:quadrotor}
\end{figure}
\Cref{fig:quadrotor} depicts the frames of reference together with the Euler angles as well as the direction of the thrust $T_j$ and angular velocity $\Omega_j$, $j \in \mathbb{I}_{1}^{4}$, of each propeller.

Based on~\cite{MaKuCo12_MAV_Dyn}, using Newton’s equations in $\mathcal{F}^i$ and Euler’s equations in $\mathcal{F}^b$, the quadrotor dynamics are
\begin{subequations}
\begin{align}
m\,\ddot{\xi} &= -m g\,\epsilon_3 + R^i_b F^b, \label{eq:quad_model_pos}\\
\dot R^i_b &= R^i_b\mathrm{S}\big(\omega_{i,b}^b\big), \label{eq:quad_model_rot}\\
J_b\,\dot{\omega}_{i,b}^b &= \tau - \omega_{i,b}^b \times \left( J_b \omega_{i,b}^b \right),
\end{align}\label{eq:quad_model}
\end{subequations}
where “$\times$” denotes the vector (cross) product, $m$ is the mass, $g$ the gravitational acceleration, $\epsilon_3=\mathrm{col}(0,0,1)$ the inertial $z$-axis, $\xi\in\mathbb{R}^3$ the position in $\mathcal{F}^i$, and $\omega_{i,b}^b=\mathrm{col}(p,q,r)$ the body angular velocity in $\mathcal{F}^b$. Let $\hat e_3=\mathrm{col}(0,0,1)$ denote the body $z$-axis; the thrust in $\mathcal{F}^b$ is $F^b=T\,\hat e_3$ with $T=\sum_{j=1}^{4}T_j$. $J_b$ is the inertia in $\mathcal{F}^b$, $\tau=\mathrm{col}(\tau_x,\tau_y,\tau_z)$ are the body torques, and $\mathrm{S}(\cdot):\mathbb{R}^3\to\mathrm{SO}(3)$ is the skew-symmetric operator.

Each rotor produces $T_j=\alpha_\mathrm{c}\Omega_j^2$ with thrust coefficient $\alpha_\mathrm{c}>0$, and drag coefficient $\beta_\mathrm{c}>0$. For the Crazyflie 2.1 ``X'' layout, let $\ell$ be the arm length and use the projected lever $d:=\ell/\sqrt{2}$. Then
\begin{equation}
\begin{bmatrix} T \\ \tau_x \\ \tau_y \\ \tau_z \end{bmatrix}
=
\begin{bmatrix}
1 & 1 & 1 & 1 \\
-d & -d &  d &  d \\
-d &  d &  d & -d \\
-\frac{\beta_\mathrm{c}}{\alpha_\mathrm{c}} & \frac{\beta_\mathrm{c}}{\alpha_\mathrm{c}} & -\frac{\beta_\mathrm{c}}{\alpha_\mathrm{c}} & \frac{\beta_\mathrm{c}}{\alpha_\mathrm{c}}
\end{bmatrix}
\begin{bmatrix} T_1 \\ T_2 \\ T_3 \\ T_4 \end{bmatrix},
\label{eq:allocation}
\end{equation}
where the signs in the last row follow rotor spin directions as in~\cref{fig:quadrotor}.

\subsection{Control architecture and prediction model}
\label{subsubsec:state_space_model}
We consider a cascaded control structure, which is the default architecture on the Crazyflie 2.1. The high-sampling-rate inner loop stabilizes attitude dynamics, while the outer loop governs position tracking. Our predictive controller operates in the outer loop and issues total thrust $T$ and body-rate $\omega_{i,b}^b$ references to the inner loop. Optimization runs off-board on a computer at 50~Hz due to the limitations of the onboard micro controller and the outer loop sends $T$ to the control mixer and $\omega_{i,b}^b$ to the inner loop, which computes $\tau$; the mixer then maps $(T,\tau)$ to rotor thrust commands $T_j,\; j \in \mathbb{I}_i^4$.

Using~\eqref{eq:quad_model_pos},~\eqref{eq:quad_model_rot}, we define a reduced-order nominal model suitable for outer-loop control. The state vector is $x=\mathrm{col}(\xi,\dot{\xi},\phi,\theta,\psi)$ and the input vector is $u=\mathrm{col}(T,p,q,r)$. 
The corresponding continuous-time state-space dynamics are
\begin{equation}
\dot{x} = f_{\mathrm{c}}\left(x, u\right) = \begin{bmatrix}
\dot{\xi} \\
 -g \epsilon_3 + \frac{1}{m} R^i_b T\hat{e}_3 \\
 p + \tan(\theta) \bigl(\sin(\phi) q + \cos(\phi) r\bigr) \\
 \cos(\phi) q - \sin(\phi) r \\
 \frac{1}{\cos(\theta)} \bigl(\sin(\phi) q + \cos(\phi) r\bigr)
\end{bmatrix}.
\label{eq:nominal_model}
\end{equation}
This baseline model $f_{\mathrm{c}}:\mathbb{R}^{n_\mathrm{x}}\times\mathbb{R}^{n_\mathrm{u}}\to\mathbb{R}^{n_\mathrm{x}}$, with $n_\mathrm{x}=9$ and $n_\mathrm{u}=4$, is discretized with a~\emph{zero-order hold} (ZOH) using~\emph{fourth-order Runge-Kutta} (RK4) at sampling time $T_{\mathrm{s}}=0.02$\,s, yielding the discrete-time predictor:
\begin{equation}
    x(k+1) = f_{\mathrm{d}}\left(x(k),u(k)\right),
    \label{eq:app:baseline-model}
\end{equation} which is used as the baseline model in the proposed predictive control scheme.

\subsection{Prediction model augmentation and residual construction}
\label{subsec:gp-training}
To capture effects not represented in the nominal model~\eqref{eq:app:baseline-model} that primarily influence the translational dynamics, only the velocity states $\dot{\xi}$ are augmented. In simulation, these effects arise mainly from the added aerodynamic drag and injected disturbances, whereas in the experiments they also reflect hardware-specific mismatches, such as the deliberate mass/inertia perturbation introduced during data collection. This leads to a residual error vector $z(k) \in \mathbb{R}^{n_{\mathrm{z}}}$ with $n_{\mathrm{z}} = 3$, where each component is modeled as an independent scalar GP $g_i \sim \mathcal{GP} \bigl( \eqref{eq:sparse_pred_mean}, \eqref{eq:sparse_pred_var} \bigr)$, for $i \in \mathbb{I}_1^{n_\mathrm{z}}$.

Let $w=\mathrm{col}(\dot{\xi},\phi,\theta,\psi,u)\in\mathbb{R}^{n_\mathrm{w}}$ be the GP input consisting of the velocity states $\dot{\xi}$, the Euler angles $\phi, \theta, \psi$, and the input $u$ excluding the position states $\xi$. At each time step $k$, the residual acceleration is computed as
\begin{equation}
z(k) = \frac{\xi_v(k+1) - \hat{\xi}_v(k+1)}{T_{\mathrm{d}}},
\label{eq:residual_error}
\end{equation}
where $\xi_v(k+1)$ and $\hat{\xi}_v(k+1)$ are the measured and predicted velocities, respectively, and $T_{\mathrm{d}}$ is the discrete sampling period. This results in an augmented model of~\eqref{eq:app:baseline-model}, expressed as
\begin{equation}
x(k+1) = f_{\mathrm{d}}\bigl(x(k), u(k)\bigr) + T_{\mathrm{d}} B_{\mathrm{z}} \hat{z}(k),
\label{eq:app_augmented_model}
\end{equation}
where the selection matrix $B_{\mathrm{z}}$ ensures the learned dynamics lie in the subspace spanned by $B_{\mathrm{z}}$, and $T_{\mathrm{d}} = T_{\mathrm{s}}$.

\subsection{Design of MPC controllers}
\label{subsec:cntl-designs}
Both controllers, LPV-MPC and NMPC, are applied with the GP-augmented prediction model~\eqref{eq:app_augmented_model} using two propagation schemes (first-order Taylor and MM, see~\cref{subsec:model-propagation}), with and without precomputed covariance propagation (\cref{subsec:precomputed-cov}). All controller configurations considered in the sequel use the same sampling time $T_{\mathrm{s}}=0.02$\,s, corresponding to a 50\,Hz outer-loop implementation in both simulation and experiments. This yields four configurations per strategy, which we denote by~\texttt{\{nl,lpv\}-\{taylor,mm\}-\{precov,cov\}} (e.g.,~\texttt{lpv-mm-precov},~\texttt{nl-taylor-cov}). In addition, we include a nominal baseline~\texttt{nl-baseline}, i.e., NMPC on the non-augmented model~\eqref{eq:app:baseline-model}.

\subsubsection{Cost function parameters}
All controllers share the same cost function~\eqref{eq:cost_func_aprox}; for~\texttt{nl-baseline} the variance trace term is omitted. The prediction horizon is fixed at $N_{\mathrm{p}}=12$, with weights $Q=\mathrm{diag}(100,100,400,40,10,10,0.1,0.1,0.1)$ and $R=\mathrm{diag}(0.1,0.1,0.1,0.1)$.

\subsubsection{Constraints}
The feasible sets $\mathcal{X}$ and $\mathcal{U}$ are imposed as:
\begin{itemize}
    \item Position: $\xi \in \mathbb{R}^3$ (unbounded),
    \item Velocity: $\dot{\xi} \in [-6.5,\; 6.5]^3$,
    \item Euler angles: $\phi, \theta, \psi \in [-70^\circ,\; 70^\circ]$,
    \item Thrust: $T \in [0.06,\; 0.64]$ according to Crazyflie limits,
    \item Angular rates: $p, q \in [-180^\circ,\; 180^\circ]$, $r \in [-20^\circ,\; 20^\circ]$.
\end{itemize}

\subsubsection{LPV-MPC iterative setup}
For~\texttt{lpv-*-*} configurations, the iterative scheme in~\cref{alg:iterated-lpv-mpc} is used with convergence tolerance $\epsilon_{\mathrm{lpv}}=0.01$ and a maximum of 12 iterations.

\subsection{Simulation study}
\label{subsec:sim_experiments}

\subsubsection{Simulation setup}
\begin{figure*}[t]
  \centering
  \includegraphics[width=0.8\linewidth]{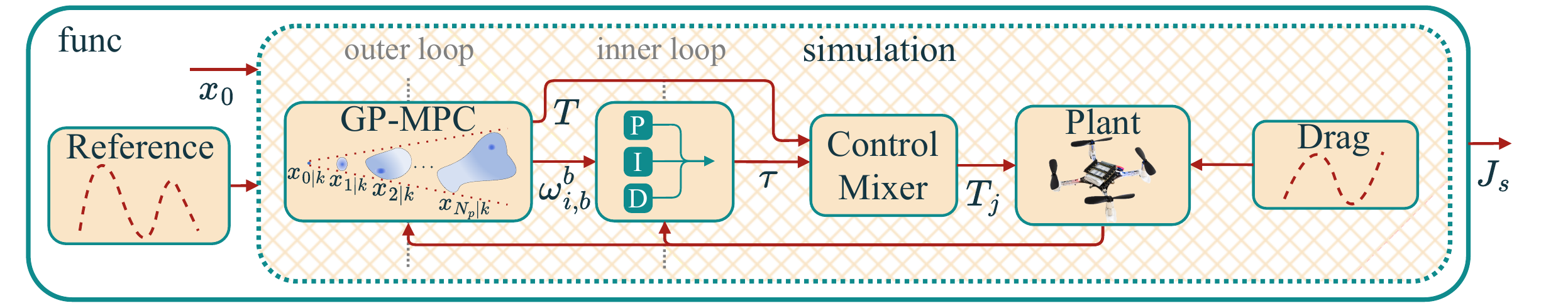}
  \caption{Schematics of the simulation environment.}
  \label{fig:simulator}
\end{figure*}

\label{subsec:sim_setup}
An open-source implementation of the proposed methodology and simulation setup is available.\footnote{\url{https://github.com/giannisbdk/gp-mpc}} All components are implemented in MATLAB; GP training (hyperparameter optimization) uses GPML~\cite{GPML}. The NMPC~\eqref{eq:nmpc} and LPV-MPC~\eqref{eq:lpv-mpc} programs are constructed in CasADi~\cite{casadi} and solved with IPOPT~\cite{IPOPT} and OSQP~\cite{OSQP}, respectively.

We simulate the Crazyflie cascaded control architecture (see~\cref{fig:simulator}), replacing the default position-velocity-attitude PID stack with MPC, while retaining the inner rate PID and the control mixer. As described in~\cref{subsubsec:state_space_model}, at 50\,Hz the outer loop sends $T$ to the firmware mixer and $\omega_{i,b}^b$ to the inner loop, which computes $\tau$; the mixer then maps $(T,\tau)$ to per-rotor thrust commands $T_j$. The full-body dynamics~\eqref{eq:quad_model} are integrated with RK4 at a numerical step size $\Delta t_{\mathrm{sim}}=0.5$\,ms to emulate continuous-time evolution between controller updates, whereas the outer-loop controllers in~\cref{subsec:cntl-designs} are updated at 50\,Hz, matching the experimental implementation in~\cref{subsec:real_experiments}.
System parameters (mass, inertia, actuator limits) and inner-loop PID gains match those in the Crazyflie firmware\footnote{\url{https://github.com/bitcraze/crazyflie-firmware}~\cite{Bitcraze}}.

To enhance fidelity, we include a simple aero-drag term computed in the body frame,
\[
F_{\mathrm{aero}}^{b} \;=\; -\Big(\textstyle\sum_{j=1}^{4}\Omega_j\Big)\,K_{\mathrm{aero}}\,R^{b}_{i}\,\dot{\xi},
\]
and use it in the translational dynamics~\eqref{eq:quad_model_pos} as $R^{i}_{b}(F^{b}+F_{\mathrm{aero}}^{b})$.

\begin{figure}[t] 
  \centering
  \includegraphics[width=0.6\linewidth]{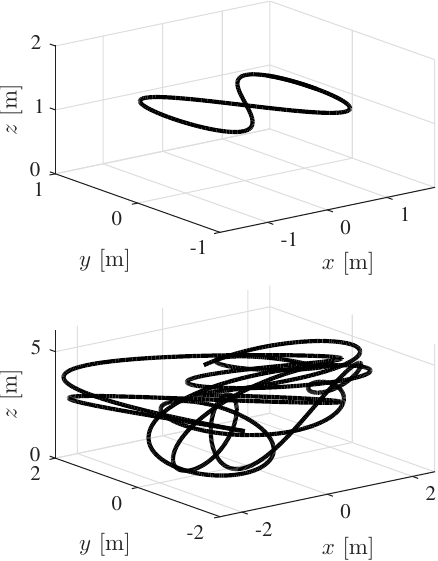}
  \caption{Lemniscate (top) and random (bottom) reference trajectories used for assessing tracking performance and for generating simulated data to train the GP, respectively.}
  \label{fig:3d_refs}
\end{figure}

Tracking is evaluated on a lemniscate reference $r(k) = \mathrm{col}\!\bigl(1.2\cos(1.3\sqrt{2}\,k),\; 1.2\sin(1.3\sqrt{2}\,k)\cos(0.77\sqrt{2}\,k),\; 1.2 + 0.02\sin(1.3\sqrt{2}\,k)\bigr)$,
as in~\cref{fig:3d_refs} (top). Performance is quantified by the root-mean-square error (RMSE)
\begin{equation}
J_{\mathrm{s}} \;=\; \sqrt{\frac{1}{N_{\mathrm{s}}}\,\sum_{k=0}^{N_{\mathrm{s}}-1} \lVert e(k) \rVert_2^2},
\label{eq:tune:cost}
\end{equation}
where $e(k) = \xi(k) - r(k) = \mathrm{col}\bigl(e_x(k), e_y(k), e_z(k)\bigr)$, over a horizon of $N_{\mathrm{s}}$ steps. Unless stated otherwise, simulations run for $10$\,s.

\subsubsection{Data collection and training}
\label{subsubsec:data_collection} 

Training data are collected by executing~\texttt{nl-baseline} on a randomized polynomial trajectory that excites the state across the operational envelope (see~\cref{fig:3d_refs}, bottom). A stochastic disturbance $v(k)\sim\mathcal N(0_3,\Sigma_v)$ with $\Sigma_v=0.1 I_3$ is injected into the translational dynamics~\eqref{eq:quad_model_pos} during data collection. Using the residual definition~\eqref{eq:residual_error}, we form $\mathcal D=\{(w_i,z_i)\}_{i=1}^N$, where $i$ denotes the sample index associated with the time-indexed quantities $(w(k),z(k))$. The construction of $(w_i,z_i)$ and the GP training procedure follow~\cref{subsec:gp-training}, yielding the augmented predictor~\eqref{eq:app_augmented_model} with GP corrections to the velocity dynamics.

\begin{figure}[tbp]
  \centering
  \includegraphics[width=1\linewidth]{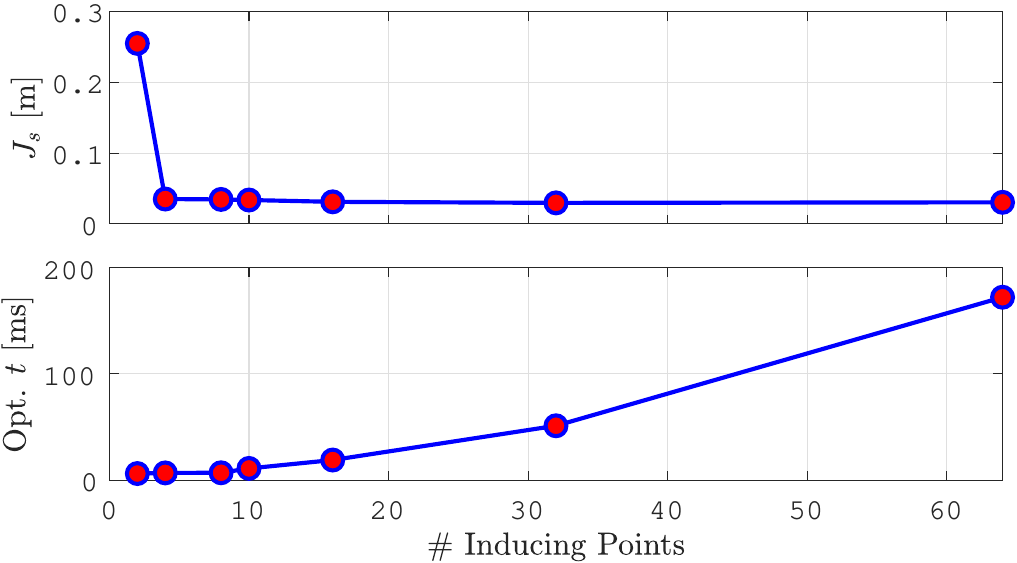}
  \caption{Simulation RMSE $J_\mathrm{s}$ by~\eqref{eq:tune:cost} (top) and average optimization time (bottom) when using the~\texttt{lpv-mm-precov} controller based on a sparse GP trained with the number of inducing points indicated at the horizontal axis.}
  \label{fig:inducing_analysis}
\end{figure}

\begin{figure}[t]
  \centering
  \includegraphics[width=1\linewidth]{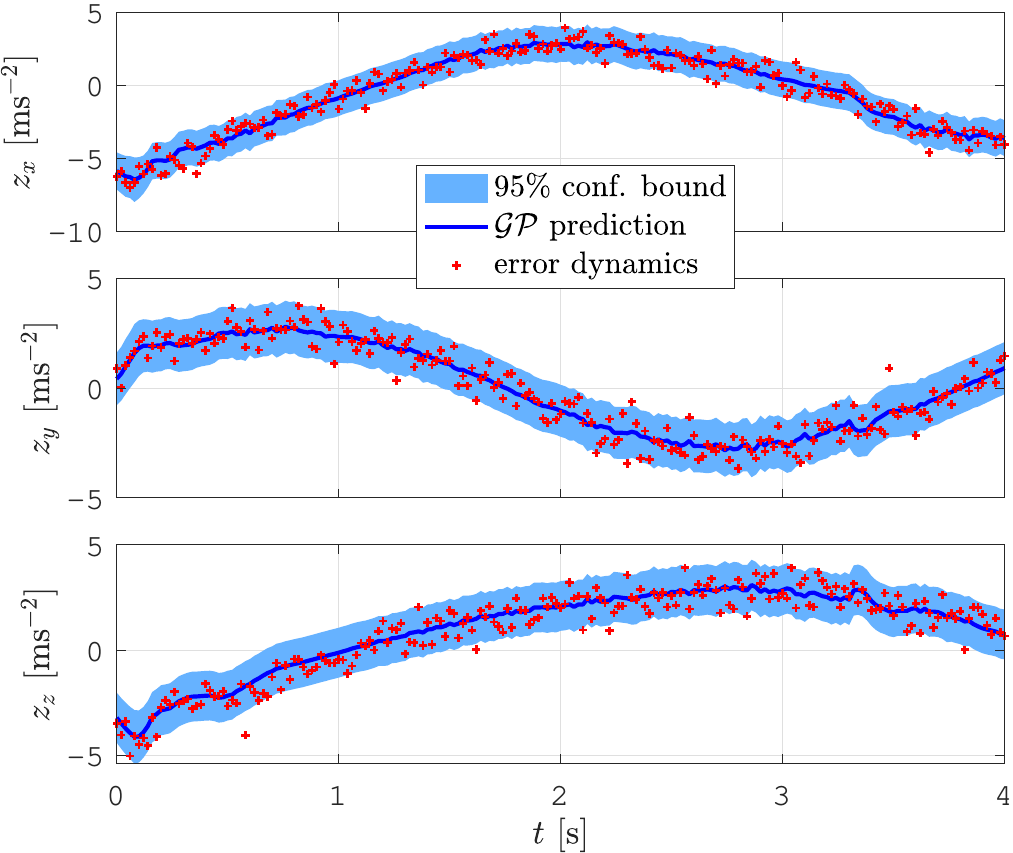}
  \caption{Sparse GP predictions along with their $95\%$ confidence bounds on a subset of the collected data.}
  \label{fig:gp_pred}
\end{figure}

To determine a suitable number of inducing points, we evaluate the trade-off between tracking accuracy (RMSE~\eqref{eq:tune:cost}) and average solver time using the~\texttt{lpv-mm-precov} controller on the lemniscate trajectory. As shown in~\cref{fig:inducing_analysis}, the RMSE saturates beyond four inducing points while solve time continues to increase; hence we use $m=4$ inducing points per GP. The learned GP predictions---mean and $95\%$ confidence bounds---on a subset of the training data are shown in~\cref{fig:gp_pred}.

\subsubsection{Simulation results}
\label{subsec:sim:performance}
\begin{figure}[t]
  \centering
  \includegraphics[width=1\linewidth]{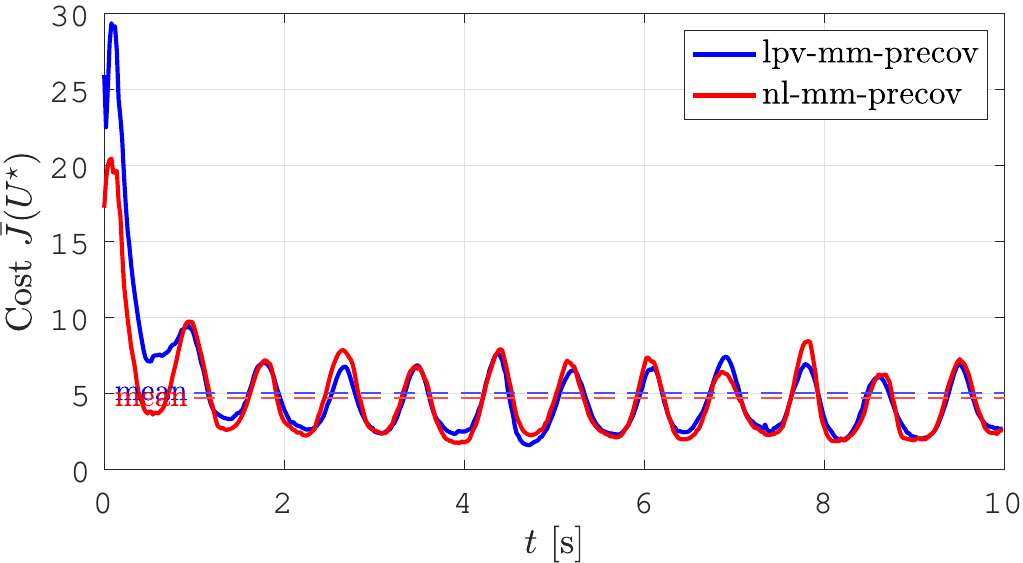}
  \caption{Overall cost for each simulated time-step $k$ between the LPV-MPC and NMPC solution using both MM and pre-computing covariance propagation. The dashed lines indicate the mean cost of each controller optimization problem.}
  \label{fig:cost_analysis}
\end{figure}

\begin{figure*}[t]
    \centering
    \includegraphics[width=1\linewidth]{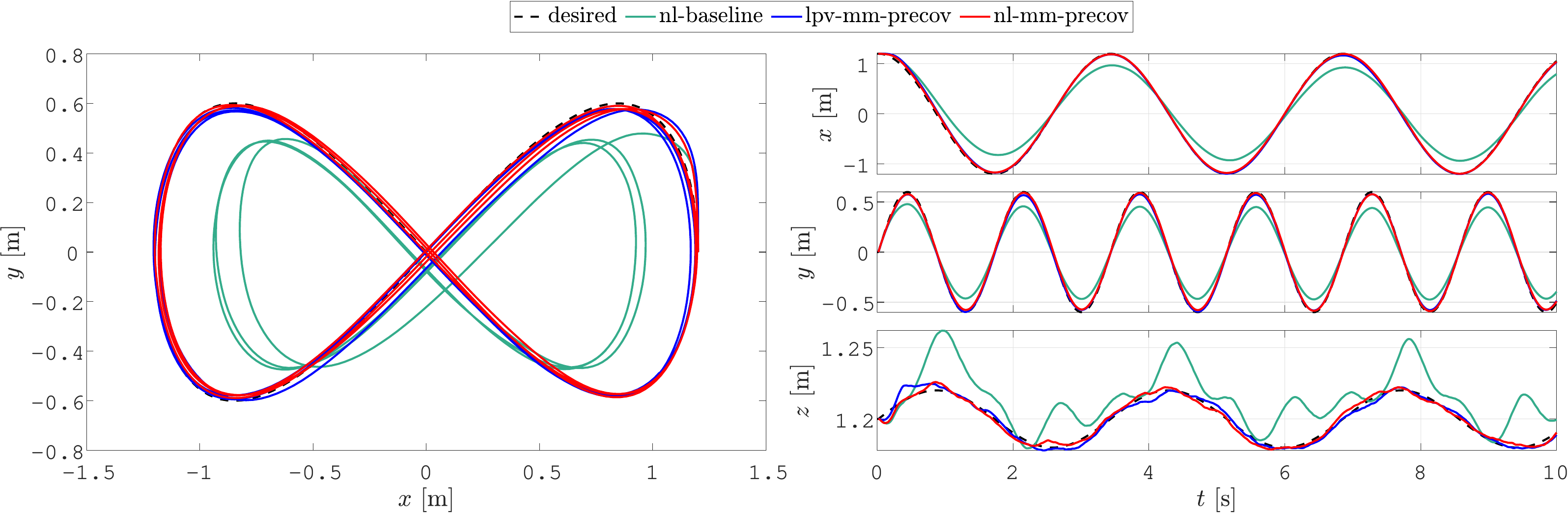}
    \caption{Closed-loop lemniscate tracking performance for the positional states $\xi$ with~\texttt{\{nl,lpv\}-mm-precov} predictive control designs.}
    \label{fig:pos_tracking}
\end{figure*}

\begin{table}[tbp]
\centering
\setlength\tabcolsep{3pt}
\centering
\caption{Comparison of closed-loop tracking with different MPC schemes in simulation.}
\begin{tabular}{lcccc}
\toprule
 \multirow{2}{*}{Controller} & \multirow{2}{*}{{$J_\mathrm{s}$~\eqref{eq:tune:cost} [\unit{\milli\metre}]}} & \multirow{2}{*}{{Avg. Time [\unit{\milli\second}]}} & \multicolumn{2}{c}{QP Solver} \\
\cmidrule(lr){4-5}
 &  &  &  Time [\unit{\milli\second}] & $\#$ Iters \\
\midrule
lpv-mm-precov       & 33    & 6.49      & 1.75      & 3     \\
lpv-taylor-precov   & 33    & 2.51      & 0.8       & 3     \\
nl-mm-precov        & 29    & 11.15     & -         & -     \\
nl-taylor-precov    & 29    & 8.65      & -         & -     \\  \cmidrule(lr){1-5}
lpv-mm-cov          & 34    & 113.58    & 49.84     & 2     \\
lpv-taylor-cov      & 35    & 17.3      & 7.18      & 2     \\
nl-mm-cov           & 28    & 538.81    & -         & -     \\
nl-taylor-cov       & 29    & 45.55     & -         & -     \\  \cmidrule(lr){1-5}
nl-baseline         & 230   & 8.5       & -         & -     \\
\bottomrule
\end{tabular}
\label{tab:comparison-controllers}
\end{table}

We evaluate LPV-MPC and NMPC in the eight configurations of~\cref{subsec:cntl-designs} on aggressive lemniscate tracking. Performance is measured by
\begin{enumerate*}[label=(\roman*)]
\item the \emph{positional} RMSE $J_{\mathrm{s}}$ of the tracking error in~\eqref{eq:tune:cost}, and
\item the \emph{average per-step} optimization (solver) time per MPC update.
\end{enumerate*}
As summarized in~\cref{tab:comparison-controllers}, LPV-MPC attains tracking accuracy comparable to NMPC across all settings; the per-step stage-cost trajectories also converge to similar values (\cref{fig:cost_analysis}), supporting LPV-MPC as an efficient surrogate.
For the LPV-MPC configurations, the reported solver time and iteration count correspond to the QP solves performed by OSQP, whereas the NMPC configurations are solved as nonlinear programs with IPOPT.

For practical deployment on MAVs (see~\cref{subsec:real_experiments}), we target a 50\,Hz outer-loop rate to
\begin{enumerate*}[label=(\roman*)]
\item maintain clear time-scale separation from the much faster inner rate loop, and
\item allocate a practical 20\,ms compute budget that covers solver time plus off-board communication overhead.
\end{enumerate*}
A configuration is deemed real-time if its average per-step optimization time remains below this 20\,ms limit with some margin. Under this criterion, all~\texttt{precov} variants satisfy the deadline; among covariance-propagating (\texttt{cov}) variants, only~\texttt{lpv-taylor-cov} sustains 50\,Hz, whereas~\texttt{lpv-mm-cov} can be viable with fewer iterations and/or with bigger convergence tolerance in~\cref{alg:iterated-lpv-mpc}.

On this trajectory, LPV-MPC converges in 2-3 iterations on average (the exact count depends on the tolerance in~\cref{subsec:cntl-designs}). While RMSEs are similar across configurations,~\texttt{lpv-mm-precov} offers the best accuracy-compute trade-off: it preserves more of the stochastic model structure via MM than Taylor-based counterparts while bounding compute through precomputed covariance. This makes it well suited for real-time control when uncertainty is moderate. In scenarios with higher uncertainty, full covariance propagation may be preferable, yet remains deployable across Taylor-based configurations. Accordingly,~\cref{fig:pos_tracking} illustrates high-fidelity tracking for~\texttt{*-mm-precov}, with both LPV-MPC and NMPC substantially outperforming the nominal~\texttt{nl-baseline}.

\subsection{Experimental study}
\label{subsec:real_experiments}

\subsubsection{Experimental setup}

\begin{figure}[t!]
  \centering
  \includegraphics[width=0.6\linewidth]{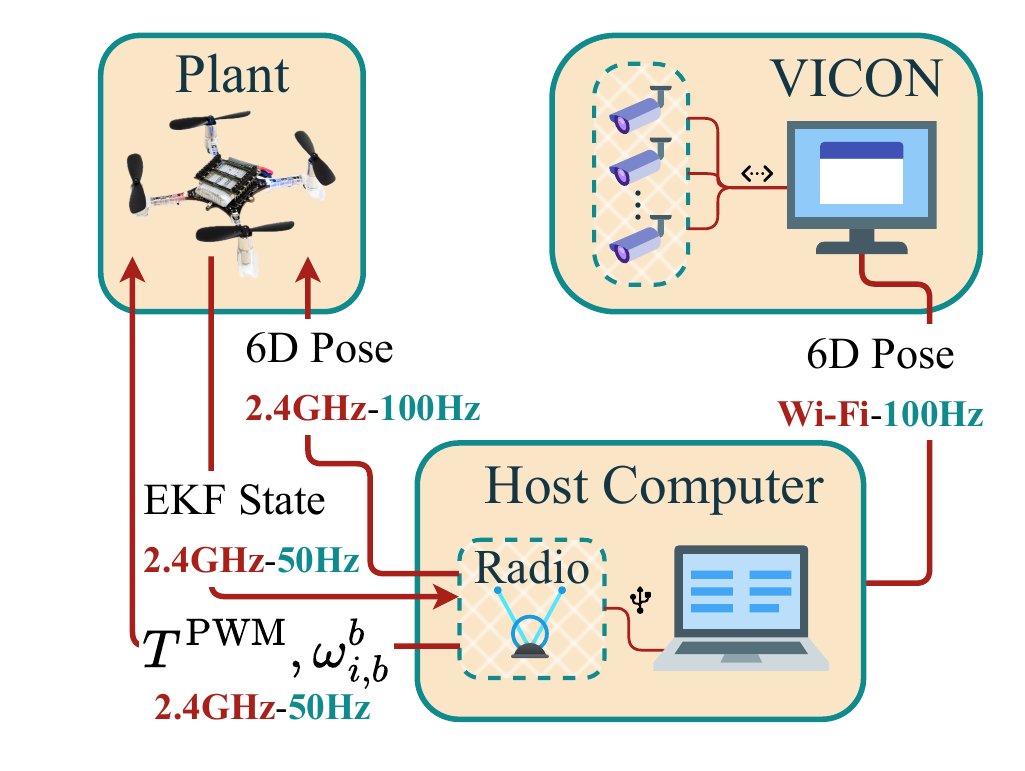}
  \caption{Schematics of the experimental environment.}
  \label{fig:exp_setup}
\end{figure}
\begin{figure}[t!]
  \centering
  \includegraphics[width=0.6\linewidth]{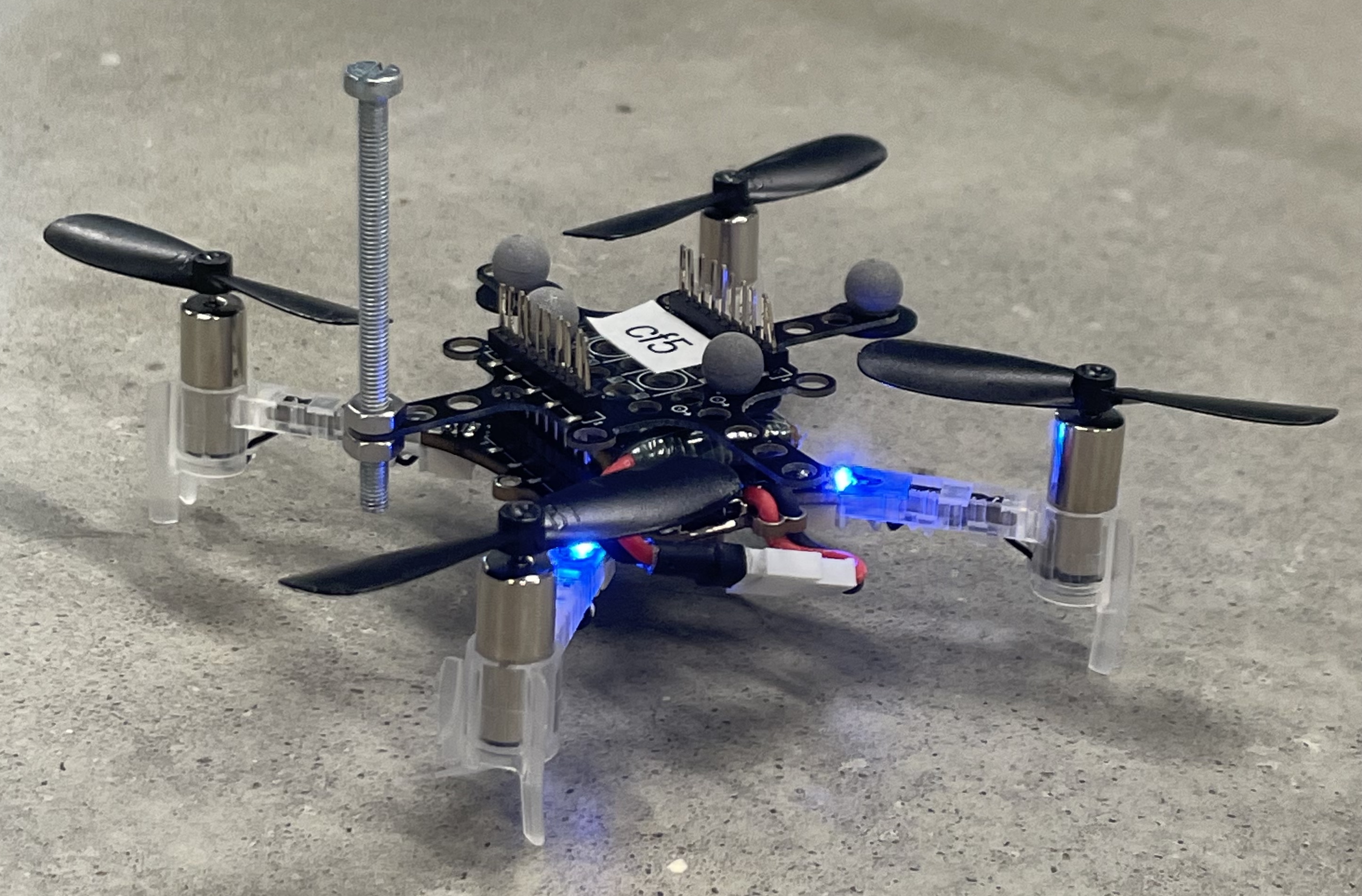}
  \caption{Crazyflie~2.1 with screw mounted on the marker deck (mass/inertia perturbation).}
  \label{fig:cf_photo}
\end{figure}

For the real-world implementation we employ a Crazyflie~2.1 tracked by a Vicon~\emph{motion-capture} (mocap) system that provides 6D pose (position $\xi$ and orientation as a quaternion $q\in\mathbb{R}^4$). A schematic of the setup and signal flow is shown in~\cref{fig:exp_setup}. The mocap system publishes pose at 100\,Hz; the host subscribes, time-stamps, and injects the pose (with respect to the inertial frame $\mathcal{F}^i$) into the on-board \emph{extended Kalman filter} (EKF) through a dedicated thread over Crazyradio. The MPC uses the estimated state $\hat{x}(k)$ of the EKF for control and sends the optimized inputs back to the Crazyflie over Crazyradio.

As in the simulation setup (\cref{subsec:sim_experiments}), the outer loop runs at 50\,Hz. At each control step $k$ the host reads $\hat{x}(k)$, solves the MPC, and transmits $(\omega_{i,b}^b,T^{\text{PWM}})$ via Crazyradio using the standard~\texttt{cflib} interface.\footnote{\texttt{send\_setpoint(roll, pitch, yawrate, thrust)}; roll/pitch are configured as~\emph{rate} setpoints in [deg/s], yawrate as rate in [deg/s], and the~\texttt{thrust} field is an integer in $[10001,60000]$. The collective thrust integer is obtained from the optimized total thrust $T^\star$ in [N] by inverting a quadratic motor thrust-PWM model under an equal rotor split ($T^\star/4$) and mapping the resulting normalized PWM to the Crazyflie integer range. We use $T_{\text{motor}}=a\,\sigma^2+b\,\sigma$ with normalized $\sigma\in[0,1]$, using the coefficients $a,b$ provided in the Crazyflie firmware; solving $a\,\sigma^2+b\,\sigma=T_{\text{motor}}$ with $T_{\text{motor}}=T^\star/4$ yields
\[
\sigma \;=\; \frac{-\,b + \sqrt{\,b^2 + 4a\,T_{\text{motor}}\,}}{2a}\in[0,1],
\]
which is then linearly mapped to the integer~\texttt{thrust} range.} The body-rate setpoint $\omega_{i,b}^b$ drive the rate inner PIDs, which produce torque-equivalent PWM terms $\tau^{\text{PWM}}$; the firmware mixer then maps the four PWM-scaled inputs $(\tau_{x}^{\text{PWM}},\tau_{y}^{\text{PWM}},\tau_{z}^{\text{PWM}},T^{\text{PWM}})$ to motor PWM $T_{j}^{\text{PWM}},\ j \in \mathbb{I}_1^4$ as
\begin{equation*}
\begin{bmatrix} T_{1}^{\text{PWM}}\\ T_{2}^{\text{PWM}}\\ T_{3}^{\text{PWM}}\\ T_{4}^{\text{PWM}} \end{bmatrix}
=
\begin{bmatrix}
1 & -1 & +1 & +1\\
1 & -1 & -1 & -1\\
1 & +1 & -1 & +1\\
1 & +1 & +1 & -1
\end{bmatrix}
\begin{bmatrix}
T^{\text{PWM}}\\[1pt] \tau_{x}^{\text{PWM}}\\[1pt] \tau_{y}^{\text{PWM}}\\[1pt] \tau_{z}^{\text{PWM}}
\end{bmatrix},
\end{equation*}
which follows the same sign pattern as the mixer in~\eqref{eq:allocation}, but it operates on PWM-scaled values. The geometric gains (projected lever arm $d$, drag $\beta_\mathrm{c}$, thrust $\alpha_\mathrm{c}$ coefficients) and motor map are absorbed into the rate inner PID controller and the thrust-PWM calibration. In that sense, the inner loop and mixer are unchanged relative to the simulated architecture.

Finally, per-step optimization time is logged on the host together with the EKF state estimate $\hat{x}(k)$ while the 20\,ms (i.e., 50\,Hz) budget includes solver time and radio overhead. The compute host is the same as in the simulation study.

\subsubsection{Data collection and training}
\begin{figure*}[t]
  \centering
  \includegraphics[width=1\linewidth]{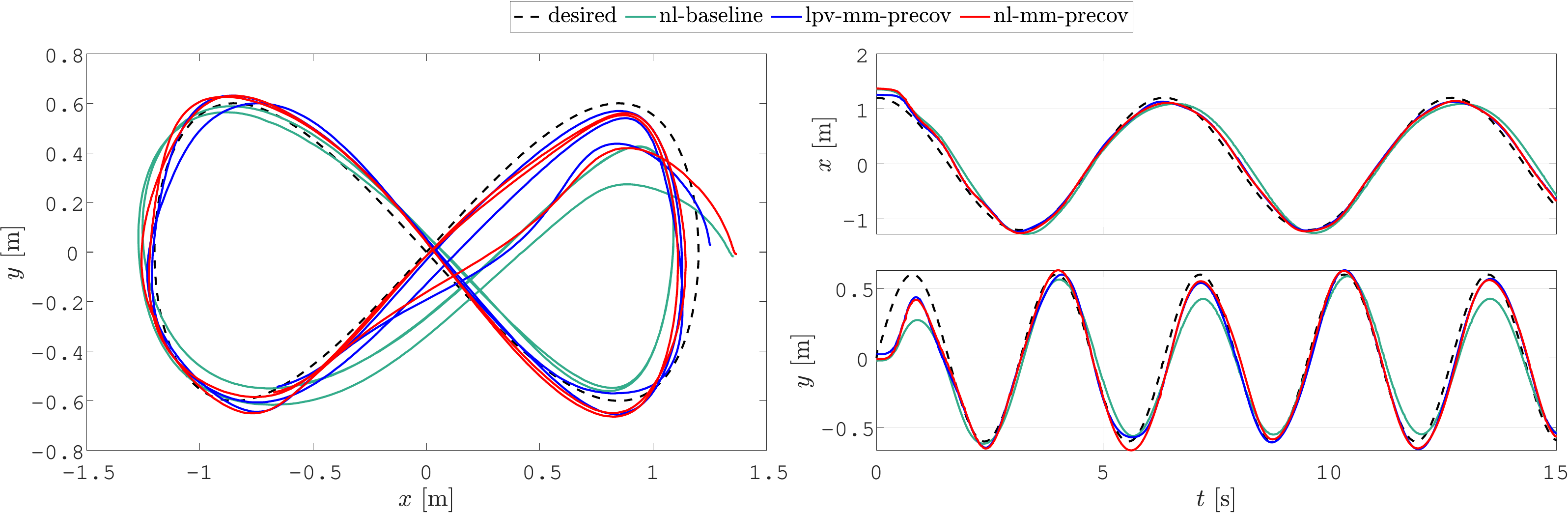}
  \caption{Real-world lemniscate tracking for the planar positional states $(\xi_x,\xi_y)$ with~\texttt{\{nl,lpv\}-mm-precov} predictive control designs.}
  \label{fig:real_tracking}
\end{figure*}

To create a deliberate model mismatch, a small screw is mounted on the mocap marker deck (\cref{fig:cf_photo}), which perturbs the mass $m$ and inertia $J_b$. As in the simulation study, GP training follows~\cref{subsec:gp-training}. The dataset is collected by executing~\texttt{nl-baseline} on the same lemniscate used later for evaluation (\cref{fig:3d_refs}, top), with the screw mounted so that the GP can learn the induced dynamics and inertia changes.

At 50\,Hz, the host computer logs the EKF state estimate $\hat x(k)$ and the commanded inputs $u(k)$. Using the residual construction in~\eqref{eq:residual_error}, we form $\mathcal D=\{(w_i,z_i)\}_{i=1}^N$ with $(w_i,z_i)$ defined as in~\cref{subsec:gp-training}; the training procedure and the resulting augmented predictor~\eqref{eq:app_augmented_model} are identical to the simulation pipeline. We fix the number of inducing points to $m=4$ per GP, as selected from the simulation trade-off.

\begin{table}[tbp]
\centering
\setlength\tabcolsep{3pt} 
\centering
\caption{Comparison of closed-loop tracking with different MPC schemes in real-world.}
\begin{tabular}{lcccc}
\toprule
 \multirow{2}{*}{Controller} & \multirow{2}{*}{{ $J_\mathrm{s}^\mathrm{xy}$~\eqref{eq:rmse_xy} [\unit{\milli\metre}]}} & \multirow{2}{*}{{Avg. Time [\unit{\milli\second}]}} & \multicolumn{2}{c}{QP Solver} \\
\cmidrule(lr){4-5}
 &  &  &  Time [\unit{\milli\second}] & $\#$ Iters \\
\midrule
lpv-mm-precov       & 85    & 2.86      & 1.3       & 2     \\
nl-mm-precov        & 87    & 5.16      & -         & -     \\ \cmidrule(lr){1-5}
nl-baseline         & 156   & 3.66      & -         & -     \\
\bottomrule
\end{tabular}
\label{tab:comparison-controllers-real}
\end{table}

\subsubsection{Experimental results}
We compare~\texttt{nl-baseline} against the favorable~\texttt{*-mm-precov} configurations selected in the simulation study, namely~\texttt{nl-mm-precov} and~\texttt{lpv-mm-precov}. The remaining variants were not retained for hardware validation, since the full-covariance NMPC configurations were not real-time feasible in simulation and therefore did not support a meaningful comparison between LPV-MPC and NMPC under the same deployment conditions.
For real-world evaluation we emphasize lateral performance and report the planar RMSE
\begin{equation}
J_\mathrm{s}^\mathrm{xy} = \sqrt{\frac{1}{N_{\mathrm{s}}}\,\sum_{k=0}^{N_{\mathrm{s}}-1}\!\bigl(e_x^2(k)+e_y^2(k)\bigr)},
\label{eq:rmse_xy}
\end{equation}
where $e_x(k)=\xi_x(k)-r_x(k)$ and $e_y(k)=\xi_y(k)-r_y(k)$, together with the average per-step solver time. Note, the simulation study reports the 3D positional RMSE $J_{\mathrm{s}}$ in~\eqref{eq:tune:cost}.

As shown in~\cref{tab:comparison-controllers-real}, both GP-augmented controllers track the lemniscate with substantially lower planar RMSE (about $44\%$ reduction) than the~\texttt{nl-baseline} and meet the 50\,Hz real-time budget.
For the LPV-MPC row, the reported solver statistics refer to OSQP, whereas the NMPC rows correspond to IPOPT.
Consistent with simulation, the LPV surrogate~\texttt{lpv-mm-precov} matches the accuracy of~\texttt{nl-mm-precov} while reducing average solve time by $\approx 44.6\%$ and converging in approximately two iterations on average. While the exact iteration count depends on the tolerances in~\cref{alg:iterated-lpv-mpc}, these results corroborate that the accelerated LPV formulation achieves real-time performance without sacrificing accuracy.

\Cref{fig:real_tracking} shows the $x,y$ time series and the $x$-$y$ path. The GP-augmented controllers follow the high-curvature segments more tightly—particularly around the outer lobes of the figure-eight—reflecting the learned mass/inertia perturbation from the screw. The \texttt{lpv-mm-precov} and \texttt{nl-mm-precov} trajectories are nearly indistinguishable in the plane and both improve clearly over~\texttt{nl-baseline}.

\section{Conclusions}
\label{sec:conclusions}
This paper proposed an accelerated LPV-MPC surrogate formulation for GP-augmented predictive models based on sparse GPs. Specifically, the GP-augmented predictor—both mean and variance dynamics—propagated over the prediction horizon under first-order Taylor and closed-form MM was recast into an exact, approximation-free LPV affine form using an FTC factorization. The resulting iterative scheme solves each MPC step as a sequence of QPs initialized by the previous solution, while keeping uncertainty coupled to the prediction dynamics, constraints, and the variance-dependent cost.

Evidence from a Crazyflie~2.1 study, including both high-fidelity simulation and real-world measurement results, shows that the LPV surrogate with MM and precomputed covariance matches the accuracy of the corresponding NMPC formulation while reducing average solve time by about $45\%$ and meeting the $20$\,ms budget. Across all controller configurations studied, both GP-augmented controller families improved tracking over the nominal baseline. Among the considered variants, the MM formulations with precomputed covariance offered the most favorable accuracy-compute trade-off, since they preserved more of the stochastic model structure than the Taylor-based counterparts while keeping the computational burden compatible with real-time execution under moderate uncertainty.

Future work includes a convergence analysis of the iterative LPV scheme and a formal connection to SQP.

\section*{Acknowledgments}
The authors thank Botond Ga\'al for assistance with the Crazyflie~2.1 hardware implementation and Michalis Galanis for helping with the figure schematics.

\appendix
\crefname{section}{Appendix}{Appendices}
\crefname{subsection}{Appendix}{Appendices}

\section*{Moment matching with Gaussian inputs: nominal model $f$ and sparse-GP augmentation}
\subsection{Mean under Gaussian input}
To find analytical expressions for $\bar{\mu}_z$, $\bar{\Sigma}_z$, and $\bar{\Sigma}_{f,z}$ evaluated under a Gaussian vector $w^\star \sim \mathcal{N} ( \mu_{w^\star}, \Sigma_{w^\star} )$, it is convenient to express~\eqref{eq:sparse_predictive_distr} as linear combinations of the kernel function $\kappa_i$~\eqref{eq:kernel} evaluated at the query input $w^\star$. The dual formulation is given as:
\begin{subequations}
    \begin{align}
        & \breve\mu_{z,i}(w^\ast) = K_{\breve w,i}^\intercal(w^\ast) \breve\alpha_i =\sum_{\tau=1}^{M} \kappa_i(w^\ast, \breve w_\tau) \breve\alpha_{i,\tau}  \label{eq:sp_pred_mean_2}, \\
        & \breve\sigma_{z,i}^2(w^\ast) = \kappa_i(w^\ast, w^\ast) - \sum_{\tau=1}^{M}\sum_{\tilde\tau=1}^{M}\sum_{\bar\tau=1}^{N} \Biggl(\Bigl([K_{\breve w \breve w,i}]^{-1}_{\tau, \tilde{\tau}} \notag \\
            & \qquad \qquad \; - [\mathcal{S}_{\breve{w}w,i}]^{-1}_{\tau, \tilde{\tau}, \bar{\tau}} \Bigr)  \kappa_i(w^\ast, \breve w_\tau) \kappa_i(w^\ast, \breve w_{\tilde\tau})\Biggr). \label{eq:sp_pred_var_2}
    \end{align}
\end{subequations}
Under this formulation for sparse GPs, the predictive mean $\bar{\mu}_z$ and the off-diagonal cross-covariances $\bar{\Sigma}_{z,i,j}$, $\forall i \neq j$, can still be derived following the same analytic procedure as in the full GP case~\cite[Section~2.3.2]{De_PILCO_PhD}, assuming modified weights $\breve \alpha$ and pseudo-inputs $\breve w$. As such, these expressions are not repeated here for brevity and are given in~\eqref{eq:gp_emm_gaussian_means} and~\eqref{eq:gp_emm_gaussian_covars}, respectively.
Additional derivation details are provided in~\cite[Appendix~A]{BadakisMSc}.
In the following, we explicitly derive the diagonal predictive variances $\bar{\Sigma}_{z,i,i}$ and the cross-covariance $\bar{\Sigma}_{f,z}$.

\subsection{Variance under Gaussian input}
Let $\bigl\{\mathrm{var}_{z,{w^{\star}}} \{ z_i | w^\star \}\bigr\}_{i=1}^{n_{\mathrm{z}}}$ denote the set of diagonal entries of the predictive covariance $\bar{\Sigma}_{z,i,i}$. By the law of total variance, for each output \(i\),
\begin{subequations}
\begin{align}
    &\mathrm{var}_{z,w^{\star}} \left \{ z_i | \mu_{w^{\star}}, \Sigma _{w^{\star}} \right\} \begin{aligned}[t] &= \mathbb{E}_{w^{\star}} \bigl\{\mathrm{var}_z\{z_i|w^{\star}\} \big| \mu_{w^{\star}}, \Sigma_{w^{\star}} \bigr\} \nonumber \\
    & + \mathrm{var}_{w^{\star}} \bigl\{\mathbb{E}_z \{z_i | w^{\star}\} \big| \mu_{w^{\star}} ,\Sigma_{w^{\star}} \bigr\} \nonumber \end{aligned} \\
    &\stackrel{\text{\eqref{eq:sparse_pred_mean},~\eqref{eq:sparse_pred_var}}}{=} \mathbb{E}_{w^{\star}} \bigl\{\breve \sigma_{z,i}^2(w^{\star})\bigr\} + \mathrm{var}_{w^{\star}} \bigl\{\breve \mu_{z,i}(w^{\star}) \bigr\} \nonumber \\
    &= \mathbb{E}_{w^{\star}} \bigl\{\breve \sigma_{z,i}^2(w^{\star})\bigr\} + \mathbb{E}_{w^{\star}} \bigl\{\breve \mu_{z,i}^2(w^{\star})\bigr\} - \mathbb{E}_{w^{\star}}^2 \bigl\{\breve \mu_{z,i}(w^{\star})\bigr\}. \nonumber
\end{align}
\end{subequations}

Since the sparse GP mean in~\eqref{eq:sp_pred_mean_2} maintains the same linear combination structure as in the full GP case (with pseudo-inputs $\breve w$ and weights $\breve \alpha$), the terms $\mathbb{E}_{w^{\star}} \bigl\{\mu_{z,i}^2(w^{\star})\}$ and $\mathbb{E}_{w^{\star}}^2 \bigl\{\mu_{z,i}(w^{\star})\}$ can be evaluated analogously, yielding $\breve \alpha_i^\intercal L_{i,i} \breve \alpha_i$ and $\bar{\mu}_{z,i}^2$ respectively, as in~\cite[Section~2.3.2]{De_PILCO_PhD}. However, the sparse predicted variance differs due to the presence of $\mathcal{S}_{\breve{w}w,i}$. The term $\mathbb{E}_{w^{\star}} \bigl\{\sigma_{z,i}^2(w^{\star})\}$ is therefore computed as:
\begin{align}
    &\mathbb{E}_{w^{\star}} \bigl\{\breve \sigma_{z,i}^2(w^{\star})\bigr\} \stackrel{\eqref{eq:sp_pred_var_2}}{=} \nonumber \\
    &=\int \Biggl[ \sigma_i^2 + \sigma_{v,i}^2  - \sum_{\tau=1}^{M}\sum_{\tilde\tau=1}^{M}\sum_{\bar\tau=1}^{N} \Biggl( \Bigl([K_{\breve w \breve w,i}]^{-1}_{\tau, \tilde{\tau}} - [\mathcal{S}_{\breve{w}w,i}]^{-1}_{\tau, \tilde{\tau}, \bar{\tau}} \Bigr)  \nonumber \\
    & \quad \cdot \kappa_i(w^\ast, \breve w_\tau) \kappa_i(w^\ast, \breve w_{\tilde\tau})\Biggr) \mathcal{N}(\mu_{w^{\star}}, \Sigma_{w^{\star}}) \Biggr] \mathrm{d}w^{\star} \nonumber \\
    &= \sigma_i^2 + \sigma_{v,i}^2 + \mathrm{Tr}\Biggl(\Bigl(K_{\breve w \breve w,i}^{-1} - \mathcal{S}_{\breve{w}w,i}^{-1} \Bigr) \sum_{\tau=1}^{M} \sum_{\tilde\tau=1}^{M} \nonumber \\
    & \quad \cdot \underbrace{\int \kappa_i(w^{\star},\breve w_\tau) \kappa_i(w^{\star}, \breve w_{\tilde\tau}) \mathcal{N}(\mu_{w^{\star}},\Sigma_{w^{\star}})\mathrm{d}w^{\star}}_{\left[ {L}_{i,i} \right]_{\tau, \tilde{\tau}}}\Biggr),\label{eq:approx_em_var_intermidiate1}
\end{align}
where \( L_{i,i} \) is given in~\eqref{eq:gp_emm_gaussian}, following the same approach as in~\cite[Section~2.3.2]{De_PILCO_PhD}. Thus, \(\mathrm{var}_{z,w^{\star}}\{ z_i \mid \mu_{w^{\star}}, \Sigma_{w^{\star}} \}\) coincides with
\eqref{eq:gp_emm_gaussian_vars}.

\subsection{Cross-covariance under Gaussian input}
The cross-covariance $\Sigma_{f,z}$ between the nominal model $f$ and the sparse GP predictor $\hat{z}$, both conditioned on the uncertain input $w^\star$, is
\begin{align}
    \Sigma_{f,z} &= \mathbb{E}_{z,f} \left\{fz^\intercal\right\} - \mathbb{E}_{f|\omega^\star}\left\{f\right\}\mathbb{E}_{z|w^\star} \left\{z^\intercal\right\} \nonumber \\
    &= \mathbb{E}_{z,f|w^{\star}} \bigl\{f z^\intercal \big| \mu_{w^{\star}}, \Sigma_{w^{\star}} \bigr\} \nonumber \\
    &\quad - \mathbb{E}_{f|w^{\star}} \bigl\{f \big| \mu_{w^{\star}}, \Sigma_{w^{\star}} \bigr\} \mathbb{E}_{z|w^{\star}} \bigl\{z^\intercal \big| \mu_{w^{\star}}, \Sigma_{w^{\star}} \bigr\} \nonumber \\
    &= \mathbb{E}_{z,w^{\star}|w^{\star}} \bigl\{\bigl(f(\mu_{w^{\star}}) + \nabla_{w^{\star}}f(\mu_{w^{\star}}) (w^{\star} - \mu_{w^{\star}})\bigr) z^\intercal \bigr\} \nonumber \\
    &\quad - \mathbb{E}_{w^{\star}} \bigl\{f(\mu_{w^{\star}}) + \nabla_{w^{\star}}f(\mu_{w^{\star}}) (w^{\star} - \mu_{w^{\star}}) \bigr\} \bar\mu^\intercal_z \nonumber \\
    &= \nabla_{w^{\star}}f(\mu_{w^{\star}}) \Bigl\{\underbrace{ \mathbb{E}_{z,w^{\star}|w^{\star}} \bigl\{w^{\star} z^\intercal \bigr\} - \mu_{w^{\star}} \bar\mu^\intercal_z}_{\bar\Sigma_{f,z}} \Bigr\}. \nonumber \label{eq:approx_em_cross_cov_intermidiate}
\end{align}
Using the first-order Taylor expansion of \(f\) at \(\mu_{w^\star}\) (see~\eqref{eq:f_taylor_gaussian}) and $\mathbb{E}_{z|w^{\star}} \{z \mid \mu_{w^\star}, \Sigma_{w^\star}\} = \mathbb{E}_{z|w^{\star}} \{\breve \mu_z(w^\star)\} = \bar \mu_z$, the per-output columns \(\bar\Sigma_{f,z,:,i}\) follow~\cite[Sec. 2.3.2]{De_PILCO_PhD}, yielding \eqref{eq:gp_emm_gaussian_crosscovars}.

\bibliographystyle{IEEEtran}
\bibliography{references}

@article{SQP1,
    title={Sequential Quadratic Programming},
    volume={4},
    DOI={10.1017/S0962492900002518},
    journal={Acta Numerica},
    author={Boggs, Paul T. and Tolle, Jon W.},
    year={1995},
    pages={1–51}
}

@inbook{INTERIOR_POINT,
author = "P{\'o}lik, Imre and Terlaky, Tam{\'a}s",
title = "Interior Point Methods for Nonlinear Optimization",
bookTitle = "Nonlinear Optimization: Lectures given at the C.I.M.E. Summer School held in Cetraro, Italy, July 1-7, 2007",
editor = "Di Pillo, Gianni and Schoen, Fabio",
year ="2010",
publisher ="Springer Berlin Heidelberg",
pages ="215--276",
}

@article{CaCaLi23_KernelMethods,
  author={Carè, Algo and Carli, Ruggero and Libera, Alberto Dalla and Romeres, Diego and Pillonetto, Gianluigi},
  journal={IEEE Control Systems Magazine}, 
  title={Kernel Methods and Gaussian Processes for System Identification and Control: A Road Map on Regularized Kernel-Based Learning for Control}, 
  year={2023},
  volume={43},
  number={5},
  pages={69-110},
  doi={10.1109/MCS.2023.3291625}
}

@article{GaAlDa23_AnnSurvey,
  author={Gianluigi Pillonetto and Aleksandr Aravkin and Daniel Gedon and Lennart Ljung and Antônio H. Ribeiro and Thomas B. Schön},
  journal={Automatica}, 
  title={Deep networks for system identification: A survey}, 
  year={2025},
  volume={171},
  number={},
  pages={111907},
  doi={https://doi.org/10.1016/j.automatica.2024.111907}
}

@inbook{RaCh05_GP,
    author = {Rasmussen, Carl Edward and Christopher K. I. Williams},
    title = {Gaussian Processes for Machine Learning},
    publisher = {MIT Press},
    year = {2005},
    pages = {171-179}
}

@inproceedings{Ti09_VFE,
    title =     {Variational Learning of Inducing Variables in Sparse Gaussian Processes},
    author =    {Titsias, Michalis},
    booktitle = {12th International Conference on Artificial Intelligence and Statistics},
    pages =     {567--574},
    year =      {2009},
    volume =    {5},
}

@article{JoCa05_UnifiedFramework_Sparse,
  author  = {Joaquin Qui{{\~n}}onero-Candela and Carl Edward Rasmussen},
  title   = {A Unifying View of Sparse Approximate Gaussian Process Regression},
  journal = {Journal of Machine Learning Research},
  year    = {2005},
  volume  = {6},
  number  = {65},
  pages   = {1939--1959}
}

@article{BackFlipping,
  author={Antal, Péter and Péni, Tamás and Tóth, Roland},
  journal={IEEE Transactions on Control Systems Technology}, 
  title={Backflipping With Miniature Quadcopters by Gaussian-Process-Based Control and Planning}, 
  year={2024},
  volume={32},
  number={1},
  pages={3-14},
  doi={10.1109/TCST.2023.3297744}}

@article{HeKaZe20_GPMPC1,
  author={Hewing, Lukas and Kabzan, Juraj and Zeilinger, Melanie N.},
  journal={IEEE Transactions on Control Systems Technology}, 
  title={Cautious Model Predictive Control Using Gaussian Process Regression}, 
  year={2020},
  volume={28},
  number={6},
  pages={2736-2743},
  doi={10.1109/TCST.2019.2949757}
}

@inproceedings{HeLiZe18_GPMPC2,
  author={Hewing, Lukas and Liniger, Alexander and Zeilinger, Melanie N.},
  booktitle={Proc. of the European Control Conference}, 
  title={Cautious NMPC with Gaussian Process Dynamics for Autonomous Miniature Race Cars}, 
  year={2018},
  volume={},
  number={},
  pages={1341-1348},
}

@article{ToKaFo21_GPMPCDrone,
  author={Torrente, Guillem and Kaufmann, Elia and Föhn, Philipp and Scaramuzza, Davide},
  journal={IEEE Robotics and Automation Letters}, 
  title={Data-Driven MPC for Quadrotors}, 
  year={2021},
  volume={6},
  number={2},
  pages={3769-3776},
  doi={10.1109/LRA.2021.3061307}
}

@article{HybridModelling,
	author = {Kessels, Bas M. and Subrahamanian Moosath, Adarsh and Fey, Rob H. B. and van de Wouw, Nathan},
	journal = {Nonlinear Dynamics},
	title = {{AI}-based state extension and augmentation for nonlinear dynamical first principles models},
	year = {2025},
}

@article{CaArWe19_MPCArm,
  author={Carron, Andrea and Arcari, Elena and Wermelinger, Martin and Hewing, Lukas and Hutter, Marco and Zeilinger, Melanie N.},
  journal={IEEE Robotics and Automation Letters}, 
  title={Data-Driven Model Predictive Control for Trajectory Tracking With a Robotic Arm}, 
  year={2019},
  volume={4},
  number={4},
  pages={3758-3765},
  doi={10.1109/LRA.2019.2929987}
}

@article{PoPeTo23_LPVGPMPC,
  author={Polcz, Péter and Péni, Tamás and Tóth, Roland},
  journal={IET Control Theory \& Applications}, 
  title={Efficient implementation of Gaussian process–based predictive control by quadratic programming}, 
  year={2023},
  volume={17},
  number={8},
  pages={968-984},
  doi={10.1109/LRA.2019.2929987}
}

@article{AmAnAn23_ZOH,
  author={Amon Lahr and Andrea Zanelli and Andrea Carron and Melanie N. Zeilinger},
  journal={European Journal of Control}, 
  title={Zero-order optimization for Gaussian process-based model predictive control}, 
  year={2023},
  volume={74},
  number={},
  pages={100862},
  doi={10.1109/LRA.2019.2929987}
}

@techreport{QuGiRa03_Uncertain_GP_EMM,
  author      = {Qui{\~n}onero-Candela, J. and Girard, A. and Rasmussen, C. E.},
  title       = {Prediction at an Uncertain Input for Gaussian Processes and Relevance Vector Machines: Application to Multiple-Step Ahead Time-Series Forecasting},
  institution = {Informatics and Mathematical Modelling, Technical University of Denmark},
  type        = {Tech. Rep.},
  number      = {IMM-TR-2003-18},
  address     = {Lyngby, Denmark},
  year        = {2003},
}

@inproceedings{DeHuHa_AnalyticMM,
    author = {Deisenroth, Marc Peter and Huber, Marco F. and Hanebeck, Uwe D.},
    title = {Analytic moment-based Gaussian process filtering},
    booktitle = {Proc. of the 26th Annual International Conference on Machine Learning},
    year = {2009},
    pages = {225–232},
}

@techreport{GiRaMu02_Uncertain_GP_Taylor,
  author      = {Girard, A. and Rasmussen, C. E. and Murray-Smith, R.},
  title       = {Gaussian Process Priors with Uncertain Inputs: Multiple-Step-Ahead Prediction},
  institution = {Department of Computing Science, University of Glasgow},
  type        = {Tech. Rep.},
  number      = {TR-2002-119},
  address     = {Glasgow, UK},
  year        = {2002},
}

@misc{LazarGuarantees,
      title={Stochastic MPC for Finite Gaussian Mixture Disturbances with Guarantees}, 
      author={Maico H. W. Engelaar and Micha P. P. Swaanen and Mircea Lazar and Sofie Haesaert},
      year={2024},
      eprint={2411.07887},
      archivePrefix={arXiv},
      primaryClass={eess.SY},
      url={https://arxiv.org/abs/2411.07887}, 
}

@phdthesis{De_PILCO_PhD,
    author = {Deisenroth, Marc Peter},
    title = {Efficient reinforcement learning using Gaussian processes},
    school = {Karlsruhe Institute of Technology , Computer Science},
    year = {2010},
    publisher = {KIT Scientific Publishing},
    series = {Karlsruhe Series on Intelligent Sensor-Actuator-Systems, Intelligent Sensor-Actuator-Systems Laboratory},
    volume = {9}
}

@phdthesis{Ko_Phd_FTCLPV,
    title = {Analysis and Control of Nonlinear Systems with Stability and Performance Guarantees: A Linear Parameter-Varying Approach},
    author = {Koelewijn, {Patrick Jan Willem}},
    year = {2023},
    publisher = {Eindhoven University of Technology},
    school = {Eindhoven University of Technology, Electrical Engineering},
}

@misc{SQP-LPVMPC,
      title={Unifying Sequential Quadratic Programming and Linear-Parameter-Varying Algorithms for Real-Time Model Predictive Control}, 
      author = {Floch, Krist{\'o}f and Lahr, Amon and T{\'o}th, Roland and Zeilinger, Melanie N.},
      year={2026},
      eprint={2511.09106},
      archivePrefix={arXiv},
      primaryClass={eess.SY},
      url={https://arxiv.org/abs/2511.09106}, 
}

@misc{Bitcraze,
    author = {},
    key = {},
    title = {{Bitcraze Crazyflie Micro Air Vehicles}},
    howpublished = {\url{https://www.bitcraze.io}},
    note = {{A}ccessed: 24/11/2025}
}

@mastersthesis{BadakisMSc,
  author = {Badakis, Giannis},
  title = {Agile Maneuvering of Micro Air Vehicles through {Gaussian} Process based Model Predictive Control},
  school = {Eindhoven University of Technology},
  year = {2024},
  month = nov,
  type = {Master's thesis},
}

@article{GPML,
  author  = {Carl Edward Rasmussen and Hannes Nickisch},
  title   = {Gaussian Processes for Machine Learning {(GPML)} Toolbox},
  journal = {Journal of Machine Learning Research},
  year    = {2010},
  volume  = {11},
  number  = {100},
  pages   = {3011--3015},
}

@article{IPOPT,
  author={A. Wächter and L. T. Biegler},
  journal={Mathematical Programming}, 
  title={On the Implementation of a Primal-Dual Interior Point Filter Line Search Algorithm for Large-Scale Nonlinear Programming}, 
  year={2006},
  volume={106},
  number={1},
  pages={25-57}
}

@article{OSQP,
  author  = {Stellato, B. and Banjac, G. and Goulart, P. and Bemporad, A. and Boyd, S.},
  title   = {{OSQP}: an operator splitting solver for quadratic programs},
  journal = {Mathematical Programming Computation},
  year    = {2020},
  volume  = {12},
  number  = {4},
  pages   = {637--672}
}

@article{BeScTo23_SUBSPACE,
  author={Gerben I. Beintema and Maarten Schoukens and Roland Tóth},
  journal={Automatica}, 
  title={Deep subspace encoders for nonlinear system identification}, 
  year={2023},
  volume={156},
  number={},
  pages={111210}
}

@article{RobotLearning,
	author = {Lee, Taeyoon and Kwon, Jaewoon and Wensing, Patrick M. and Park, Frank C.},
	journal = {Annual Review of Control, Robotics, and Autonomous Systems},
	pages = {311-334},
	title = {Robot Model Identification and Learning: A Modern Perspective},
	volume = {7},
	year = {2024},
}

@article{SaKaAr23_AgileMavsANN_MPC,
  author={Salzmann, Tim and Kaufmann, Elia and Arrizabalaga, Jon and Pavone, Marco and Scaramuzza, Davide and Ryll, Markus},
  journal={IEEE Robotics and Automation Letters}, 
  title={Real-Time Neural MPC: Deep Learning Model Predictive Control for Quadrotors and Agile Robotic Platforms}, 
  year={2023},
  volume={8},
  number={4},
  pages={2397-2404},
  doi={10.1109/LRA.2023.3246839}
}

@article{MaKuCo12_MAV_Dyn,
  author={Mahony, Robert and Kumar, Vijay and Corke, Peter},
  journal={IEEE Robotics \& Automation Magazine}, 
  title={Multirotor Aerial Vehicles: Modeling, Estimation, and Control of Quadrotor}, 
  year={2012},
  volume={19},
  number={3},
  pages={20-32},
  doi={10.1109/MRA.2012.2206474}
}

@article{casadi,
  author = {Joel A E Andersson and Joris Gillis and Greg Horn
            and James B Rawlings and Moritz Diehl},
  title = {{CasADi} -- {A} software framework for nonlinear optimization
           and optimal control},
  journal = {Mathematical Programming Computation},
  volume = {11},
  number = {1},
  pages = {1--36},
  year = {2019},
  publisher = {Springer},
  doi = {10.1007/s12532-018-0139-4}
}

@inproceedings{LPV-FTC-new,
    author = {E. Javier Olucha and Patrick J. W. Koelewijn and Amritam Das and Roland Tóth},
    title = {Automated Linear Parameter-Varying Modeling of Nonlinear Systems: A Global Embedding Approach},
    booktitle = {Proc. of the 6th IFAC Workshop on Linear Parameter Varying Systems},
    year = {2025}
}

\end{document}